\pdfoutput=1
\PassOptionsToPackage{dvipsnames}{xcolor}
\documentclass[final,onefignum,onetabnum]{siamart190516}
\usepackage{lipsum}
\usepackage{amsfonts}

\usepackage{graphicx}
\graphicspath{{figs/}}
\usepackage{amsmath, amssymb, mathrsfs}

\usepackage[T1]{fontenc}
\usepackage[utf8]{inputenc}
\usepackage{microtype}
\hypersetup{
  colorlinks=true,
  urlcolor=NavyBlue,
  citecolor=NavyBlue,
  linkcolor=NavyBlue,
  frenchlinks=false,
  pdfborder={0 0 0},
  naturalnames=false,
  hypertexnames=false,
  breaklinks
}
\usepackage[all]{hypcap}
\usepackage{orcidlink}

\usepackage{fontawesome}
\definecolor{twitterblue}{RGB}{64,153,255}
\newcommand{\twitter}[1]{\href{https://twitter.com/#1}{\textcolor{twitterblue}{\faTwitter}\,\tt \textcolor{twitterblue}{@#1}}}

\newsiamremark{remark}{Remark}
\newsiamremark{hypothesis}{Hypothesis}
\crefname{hypothesis}{Hypothesis}{Hypotheses}
\newsiamthm{claim}{Claim}

\newcommand{\ourtitle}{Surprises in a classic boundary-layer problem}
\newcommand{\shortauthors}{Clark, Gomes, Rodriguez-Gonzalez, Stein, Strogatz}

\headers{\ourtitle}{\shortauthors}

\title{\ourtitle\thanks{%
    Submitted to the editors July 24, 2021.
    \funding{%
    This work was supported by the National Science Foundation (DMS--1645643).}}}

\author{William~A.~Clark\thanks{%
   {\large \orcidlink{0000-0002-0817-0309}}
    Department of Mathematics, Cornell University, Ithaca, New York 14853, USA
    (\email{wac76@cornell.edu})
  }%
  \and
  Mario~W.~Gomes\thanks{%
    {\large \orcidlink{0000-0002-7593-1163}}
    Department of Mechanical Engineering, Rochester Institute of Technology, Rochester, NY 14623, USA
    (\email{mwgeme@rit.edu})
  }%
  \and
  Arnaldo~Rodriguez-Gonzalez\thanks{%
    {\large \orcidlink{0000-0002-9200-4505}}
    Sibley School of Mechanical and Aerospace Engineering, Cornell
    University, Ithaca, New York 14853, USA
    (\email{ajr295@cornell.edu},
    \twitter{Arnaldo\textunderscore AGITF})
  }%
  \and
  Leo~C.~Stein\thanks{%
    {\large \orcidlink{0000-0001-7559-9597}}
    Department of Physics and Astronomy, The University of Mississippi, University, MS 38677, USA
    (\email{lcstein@olemiss.edu},
     \twitter{duetosymmetry})
  }%
  \and
  Steven~H.~Strogatz\thanks{%
    {\large \orcidlink{0000-0003-2923-3118}}
    Department of Mathematics, Cornell University, Ithaca, New York 14853, USA
    (\email{shs7@cornell.edu},
    \href{https://stevenstrogatz.com/}{stevenstrogatz.com})
  }
}

\usepackage{amsopn}

\DeclareMathOperator{\sech}{sech}
\DeclareMathOperator{\Li}{Li}

\newcommand{\zc}{z_{c}}

\hypersetup{
  pdftitle={\ourtitle},
  pdfauthor={\shortauthors}
}

\newcommand{\figExampleSols}{%
  \begin{figure}
    \centering
    \includegraphics[width=5.in]{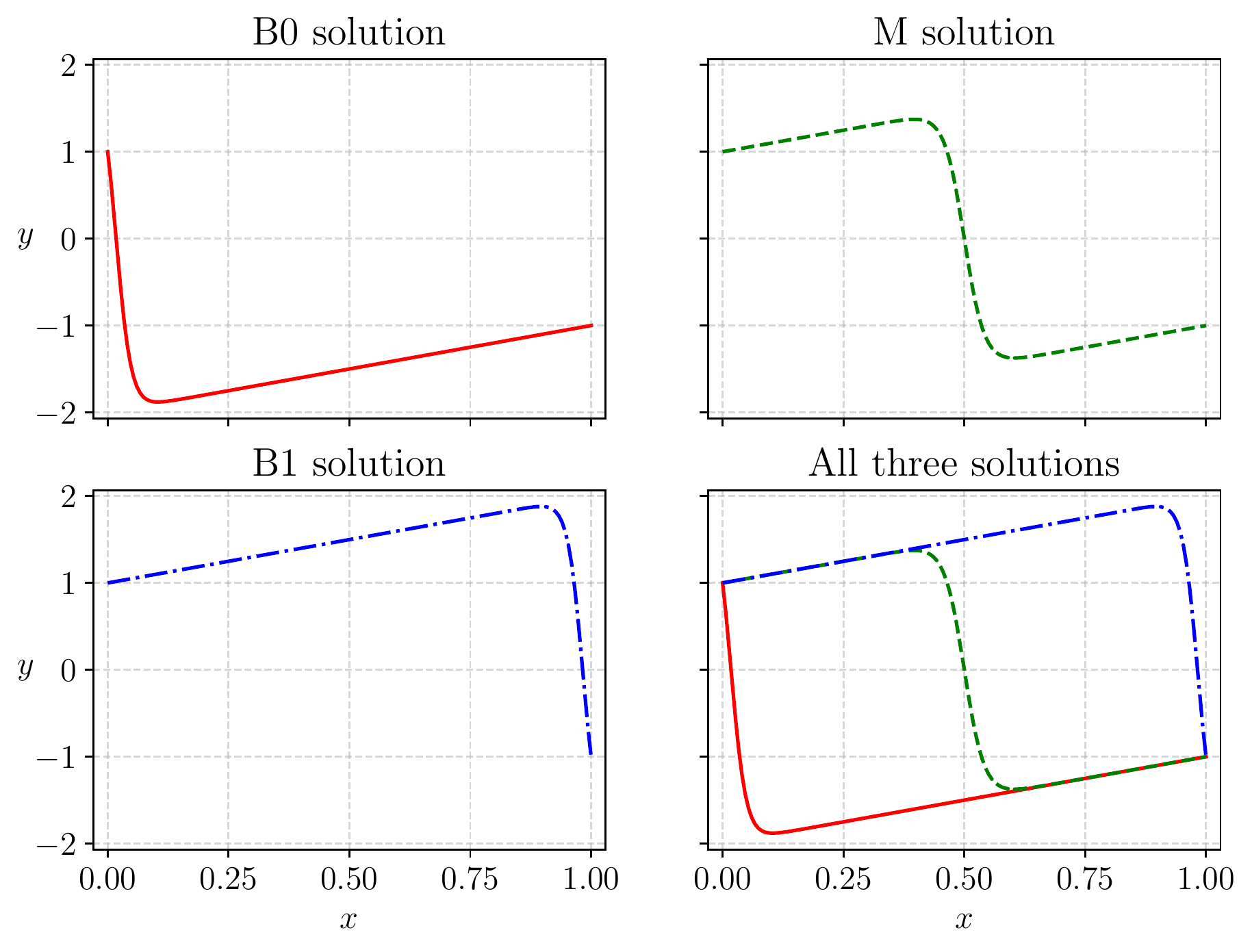}
    \caption{%
      Three numerical solutions of problem~\eqref{eq:ODE_orig} for $\epsilon=0.03$: two with
      boundary layers, and one with an interior layer.  Under the
       symmetry $(x, y) \to (1-x, -y)$, the left and right
      boundary-layer solutions B0 and B1 trade places, and the interior-layer
      solution M is unchanged.
    }
    \label{fig:example_sols_and_phase_plane}
  \end{figure}%
}

\newcommand{\figPhasePlaneSimple}{%
  \begin{figure}
    \centering
    \includegraphics[height=2.17in]{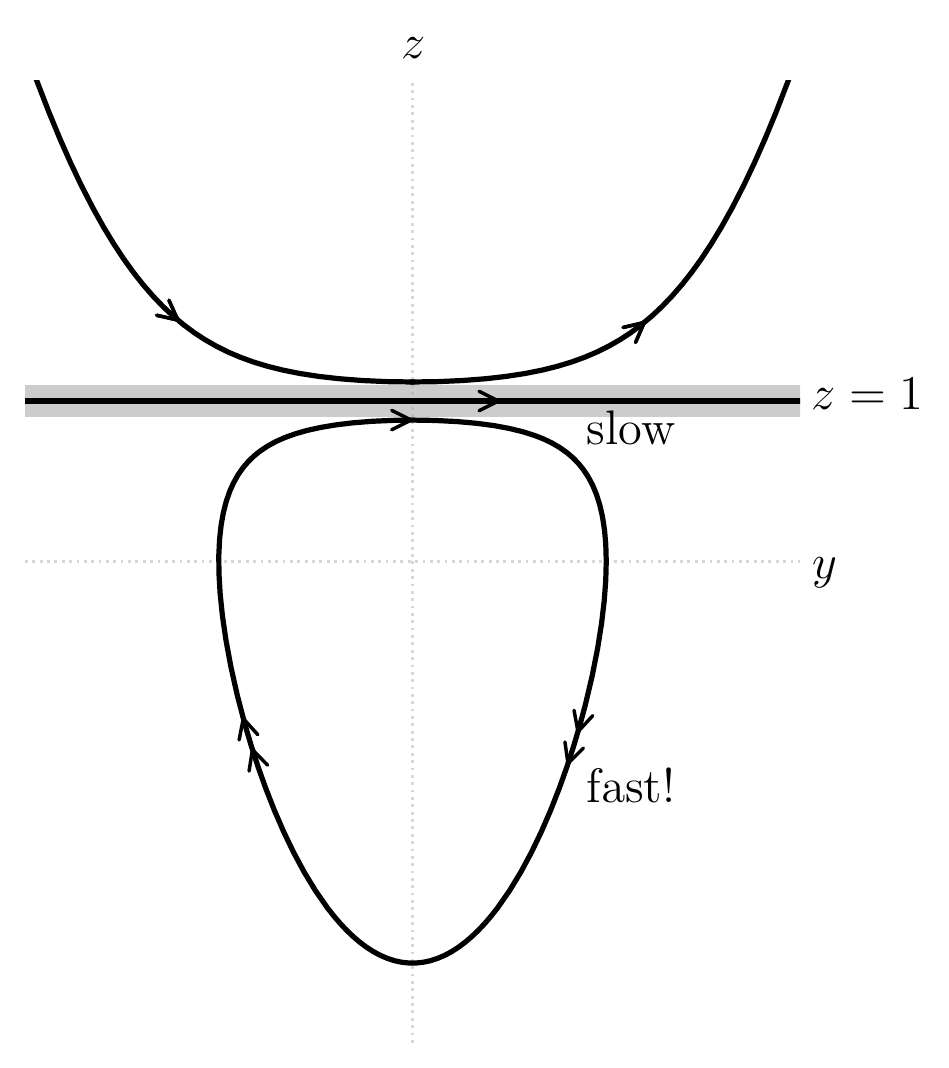}
    \includegraphics[height=2.17in]{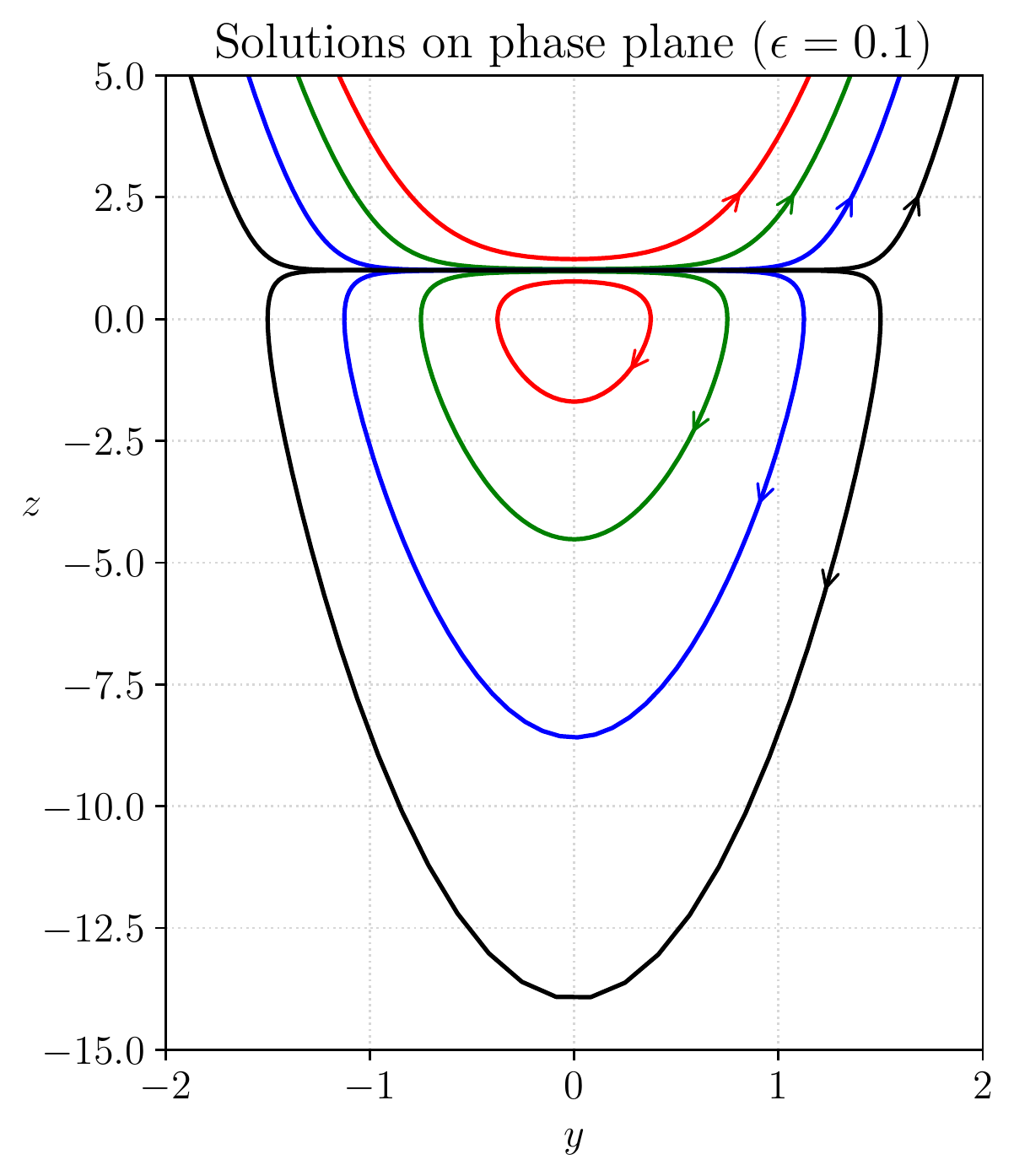}
    \caption{%
      Phase portraits for \eqref{eq:ODE_reduced} with   $\epsilon=0.1$.
      The left panel shows three qualitatively different kinds of trajectories, as well as fast and slow regions in the flow. An imaginary particle moves with $O(1)$ speed in the slow region (the thin gray strip around the invariant line $z=1$). Outside the strip, the particle  zips around much faster, with enormous vertical speeds of order $O(1/\epsilon)$.  The right panel shows a quantitatively accurate phase portrait. Notice how tightly packed the trajectories become as they squeeze into and out of the slow region near $z=1$.%
    }
    \label{fig:phase_plane_simple}
  \end{figure}%
}

\newcommand{\figPhasePlaneThreeSols}{%
  \begin{figure}
    \centering
    \includegraphics[height=2.17in]{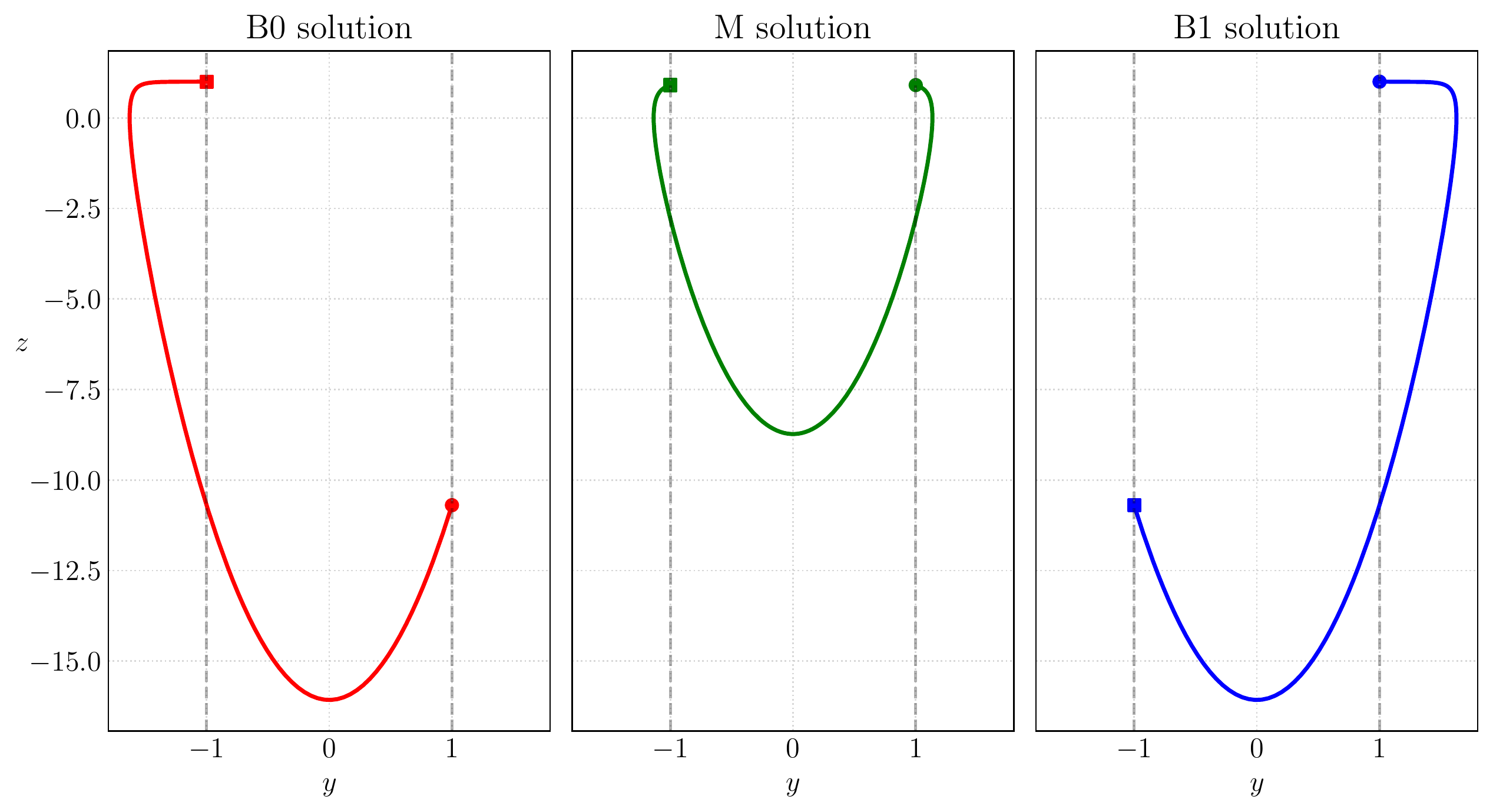}
    \caption{%
      Phase plane plots of the three solutions of the boundary-value problem~\eqref{eq:ODE_orig} for $\epsilon = 0.1$. The vertical coordinate $z$ denotes $y'$. The red, green, and blue, trajectories correspond to the B0, M, and B1 solutions having
      layers at $x=0$, $x=1/2$,and $x=1$
      respectively. The beginning ($x=0$) of each trajectory is denoted
      with a circle, and the end ($x=1$) is denoted with a square; the flow is clockwise on each trajectory. All trajectories start on the dashed vertical line $y=1$ and end on the dashed vertical line $y=-1$, corresponding to the original boundary conditions. Note that the B0 and B1 trajectories form a symmetric pair. As such, both of them lie on the same periodic orbit in the phase plane. The slow part of each trajectory (the ``outer solution'') occurs at the top, close to $z=1$, while the much faster part (the ``inner solution'' in the layer) occurs everywhere else below that.   
    }
    \label{fig:phase_plane_three_sols}
  \end{figure}%
}

\newcommand{\figUniformError}{%
  \begin{figure}
    \centering
    \includegraphics[height=3in]{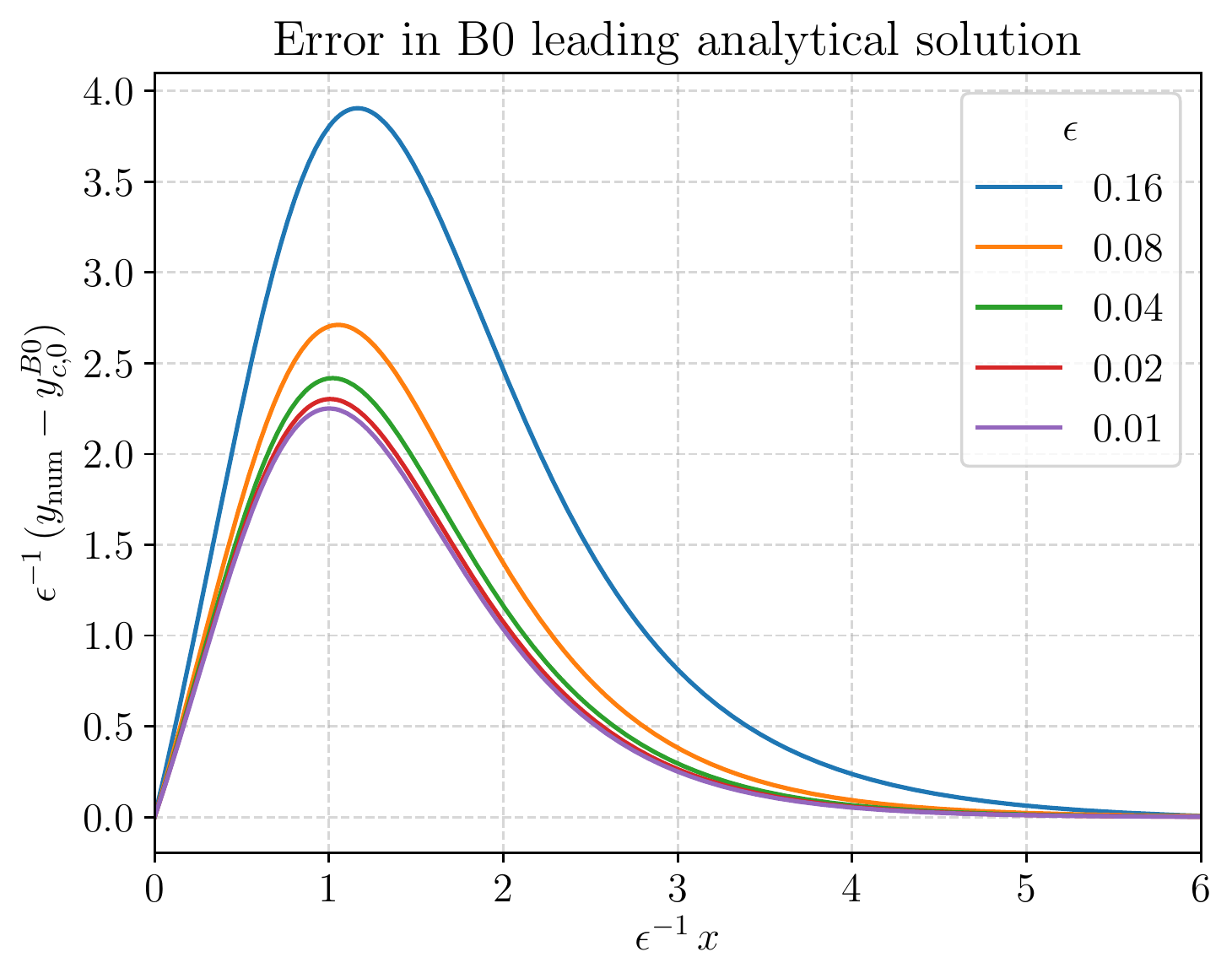}
    \caption{%
      Scaled error in the composite zeroth-order asymptotic solution \eqref{eq:yB0c0} for the case with a boundary layer at $x=0$.  The error is defined as the difference between the asymptotic solution \eqref{eq:yB0c0} and a very careful numerical solution (taken as a surrogate for the unknown exact solution). The scaled error is defined here as the error divided by $\epsilon$; this scaling is appropriate because we expect to incur errors of size $O(\epsilon)$ in a leading-order solution.  We examine the error within the boundary layer,
      $x=O(\epsilon)$, i.e.,\ the layer has a width linear in
      $\epsilon$ (the error outside of the layer is transcendentally
      small).  Notice that the composite solution is uniformly
      valid, and the error is proportional to $\epsilon$, as
      expected from a zeroth-order solution. The error behavior is similar for the other two asymptotic 
      solutions with layers at $x=1/2$ and $x = 1.$  %
    }
    \label{fig:uniform_conv}
  \end{figure}%
}

\newcommand{\figUniformErrorHigherOrder}{%
  \begin{figure}
    \centering
    \includegraphics[height=3in]{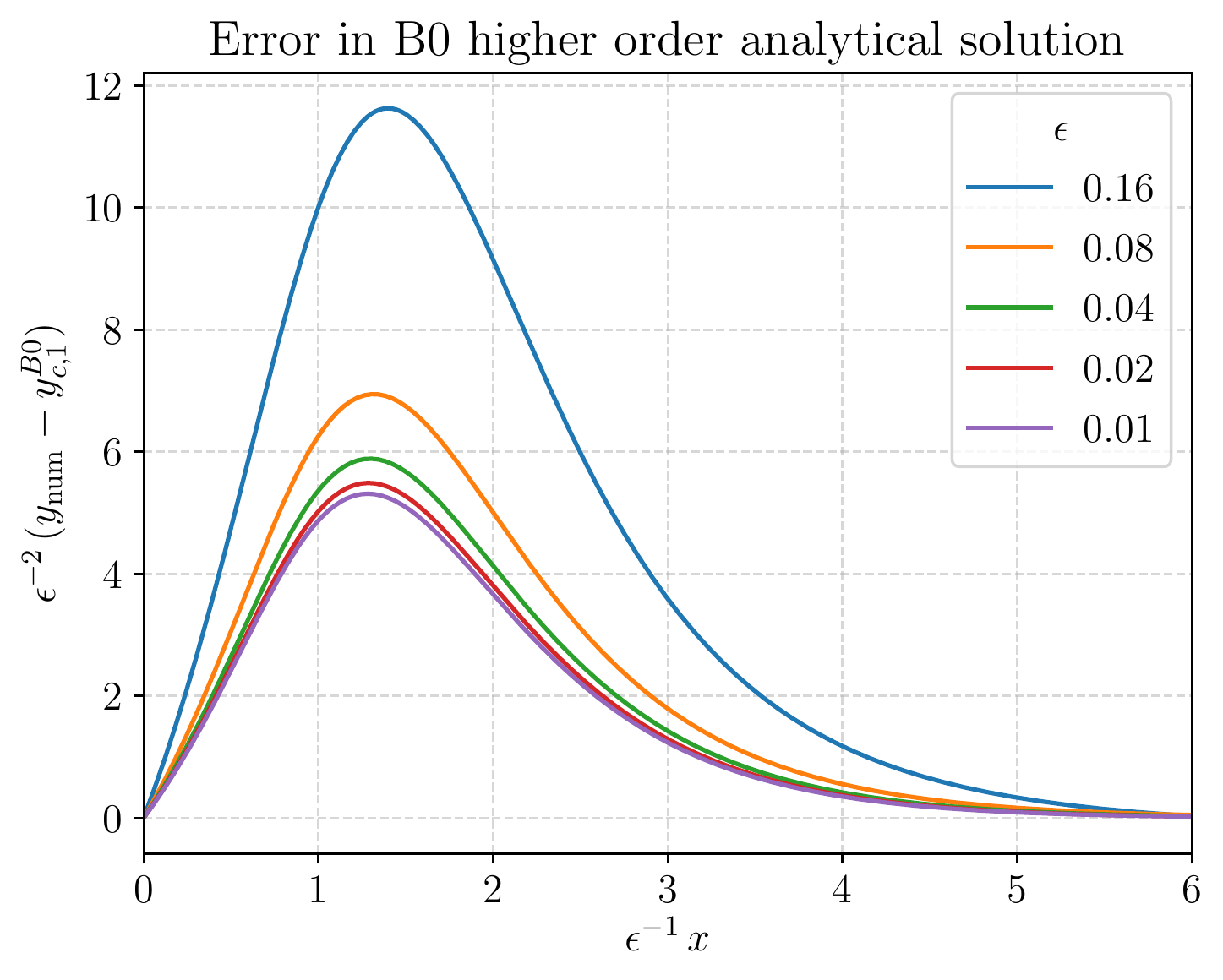}
    \caption{%
      Scaled error, now defined as error divided by $\epsilon^2$, in the composite first-order asymptotic solution for the
      case with a boundary layer at $x=0$, given in
      Appendix~\ref{sec:high-order-asympt}.
      Notice that the composite solution is uniformly
      valid, and the error is now proportional to $\epsilon^2$, as
      expected for an asymptotic solution matched to first order.%
    }
    \label{fig:uniform_conv_higher}
  \end{figure}%
}
\newcommand{\figShootingTrials}{%
  \begin{figure}
    \centering
    \includegraphics[height=2.1in]{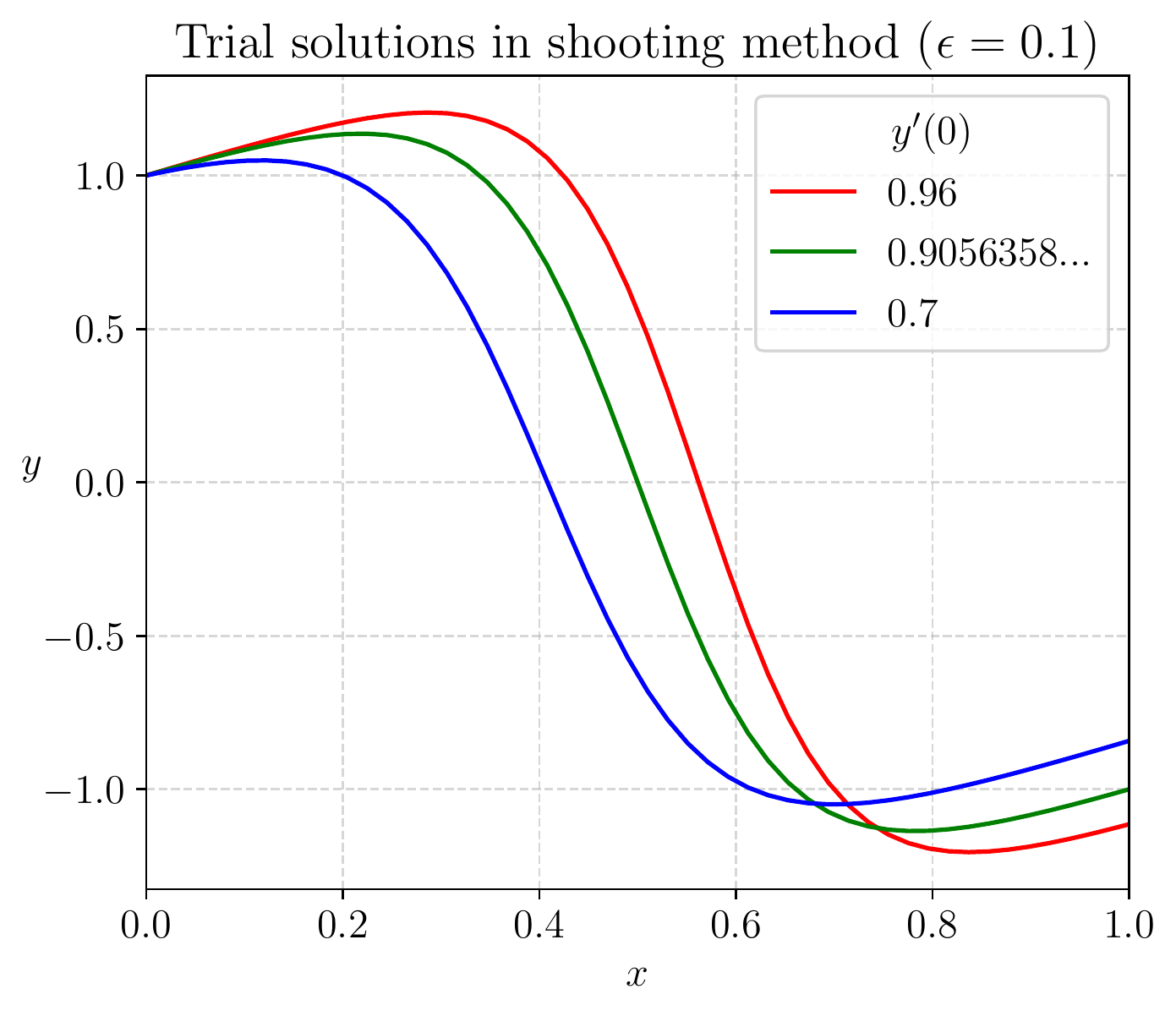}
    \includegraphics[height=2.1in]{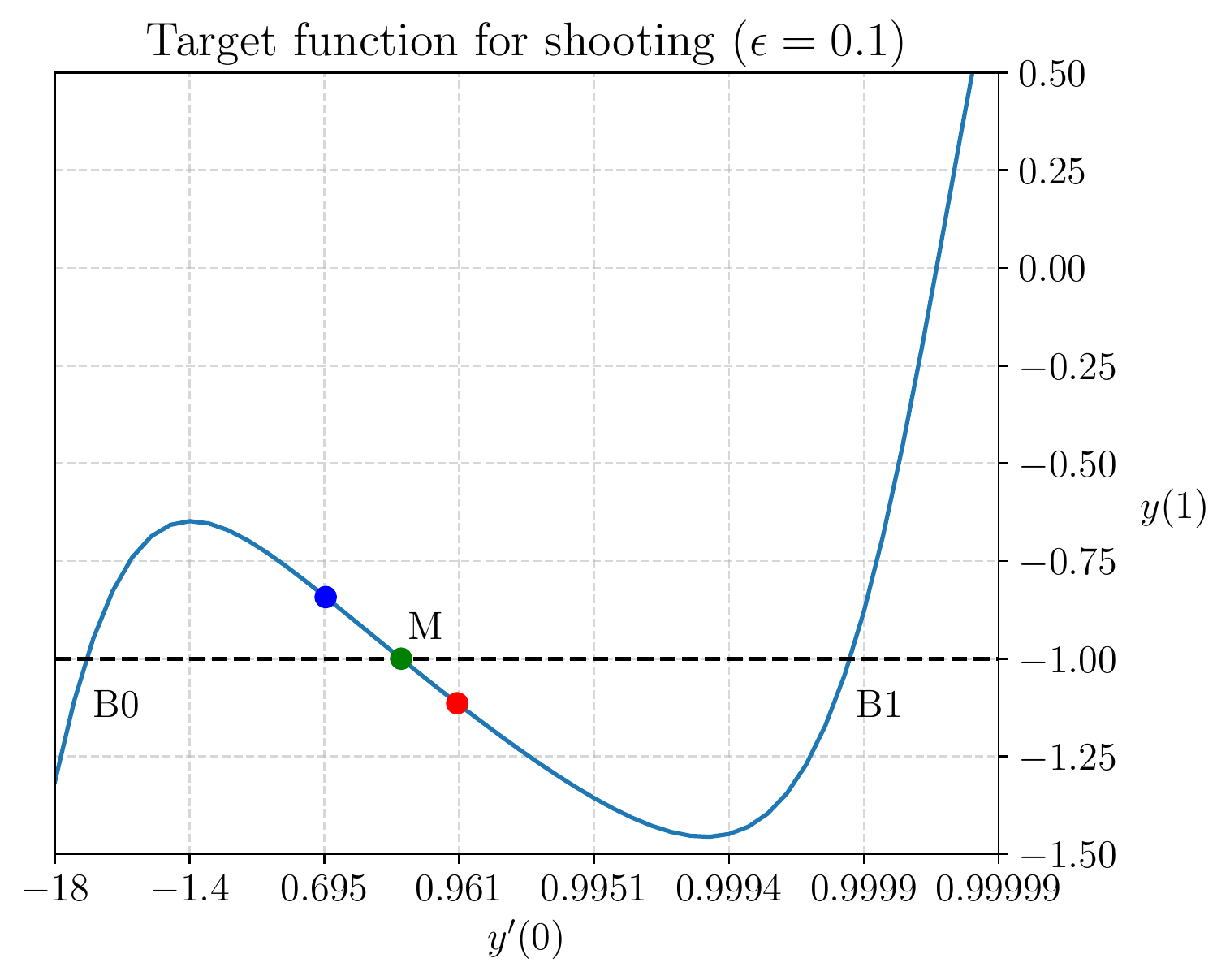}
    \caption{Solving the boundary-value problem via shooting.  Different initial
      slopes lead to different final values of $y$.
      The red, green, and blue trial solutions in the left panel
      have certain initial slopes $y'(0)$ and final $y(1)$ which
      correspond to the three filled circles plotted in the right panel.
      If we were to vary the initial slope continuously over a broader
      range, we would build up the curve in the right panel.
      Note that the horizontal axis in the right panel is uniform in
      $\log(1-y')$, to expand the exponentially bunched region near $y'(0) = 1$. The initial slopes for the B0, M, and B1 solutions are marked with letters in the right panel; they occur at the three values of $y'(0)$ where the graph of the target function crosses the dashed horizontal line $y(1)=-1$. 
    }
    \label{fig:shooting_trials}
  \end{figure}%
}

\newcommand{\figInitialSlopes}{%
  \begin{figure}
    \centering
    \includegraphics[height=3in]{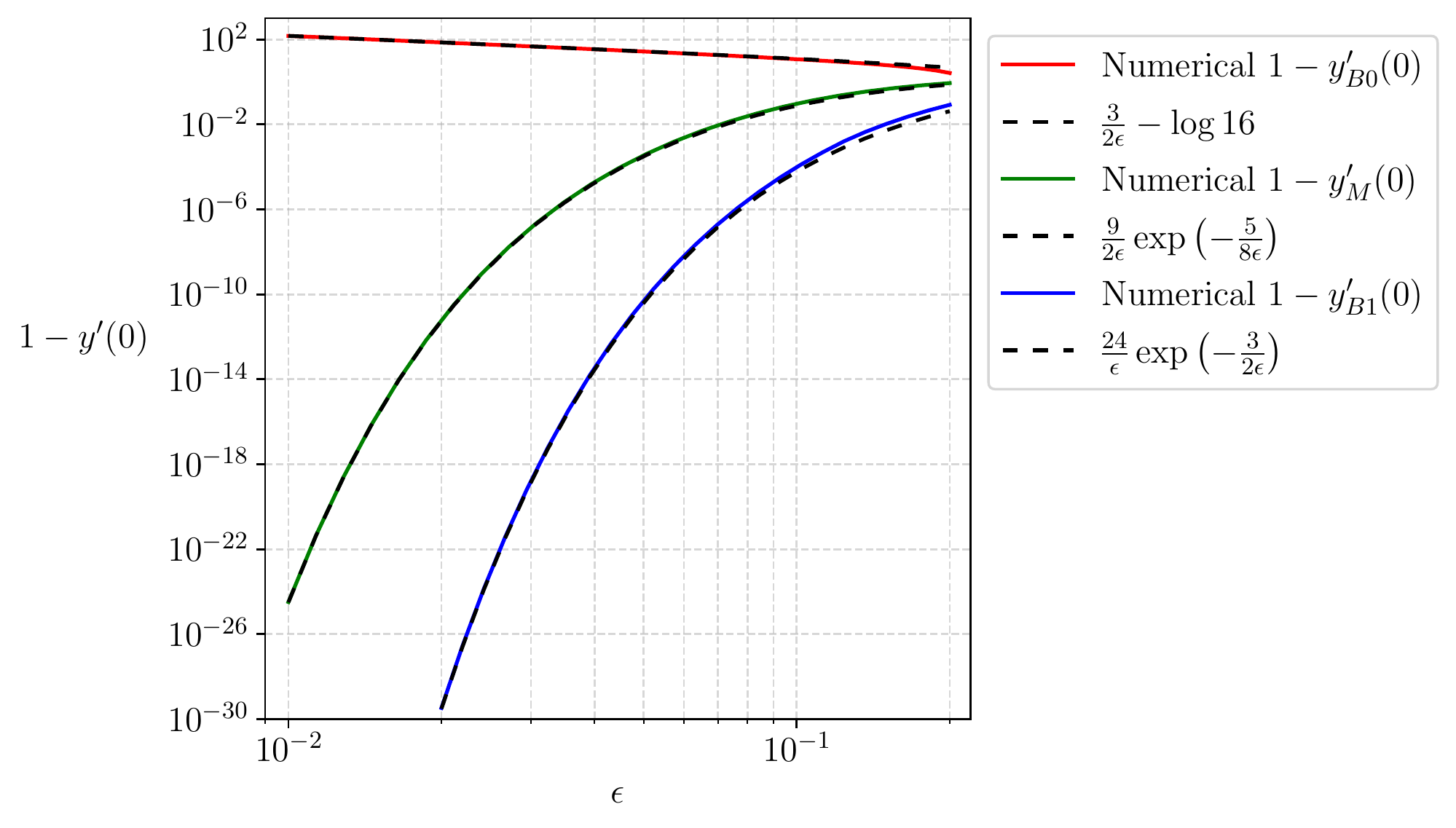}
    \caption{%
      Comparison of the initial slopes $y'(0)$ determined from the
      asymptotic solutions (dashed curves) and numerical solutions
      found via the shooting method (solid curves). The vertical axis plots $1-y'(0)$ to highlight the tiny (transcendentally small) deviation of the initial slope from 1 for two of the solutions. From top to
      bottom, the curves give the behavior of the solutions labeled B0
      (boundary layer at $x=0$), M (interior layer at $x=1/2$), and B1
      (boundary layer at $x=1$).  Notice the vast difference in
      scale for the three initial slopes.%
    }
    \label{fig:initial_slopes}
  \end{figure}%
}

\newcommand{\figPitchfork}{%
  \begin{figure}[t]
    \centering
    \includegraphics[width=4.in]{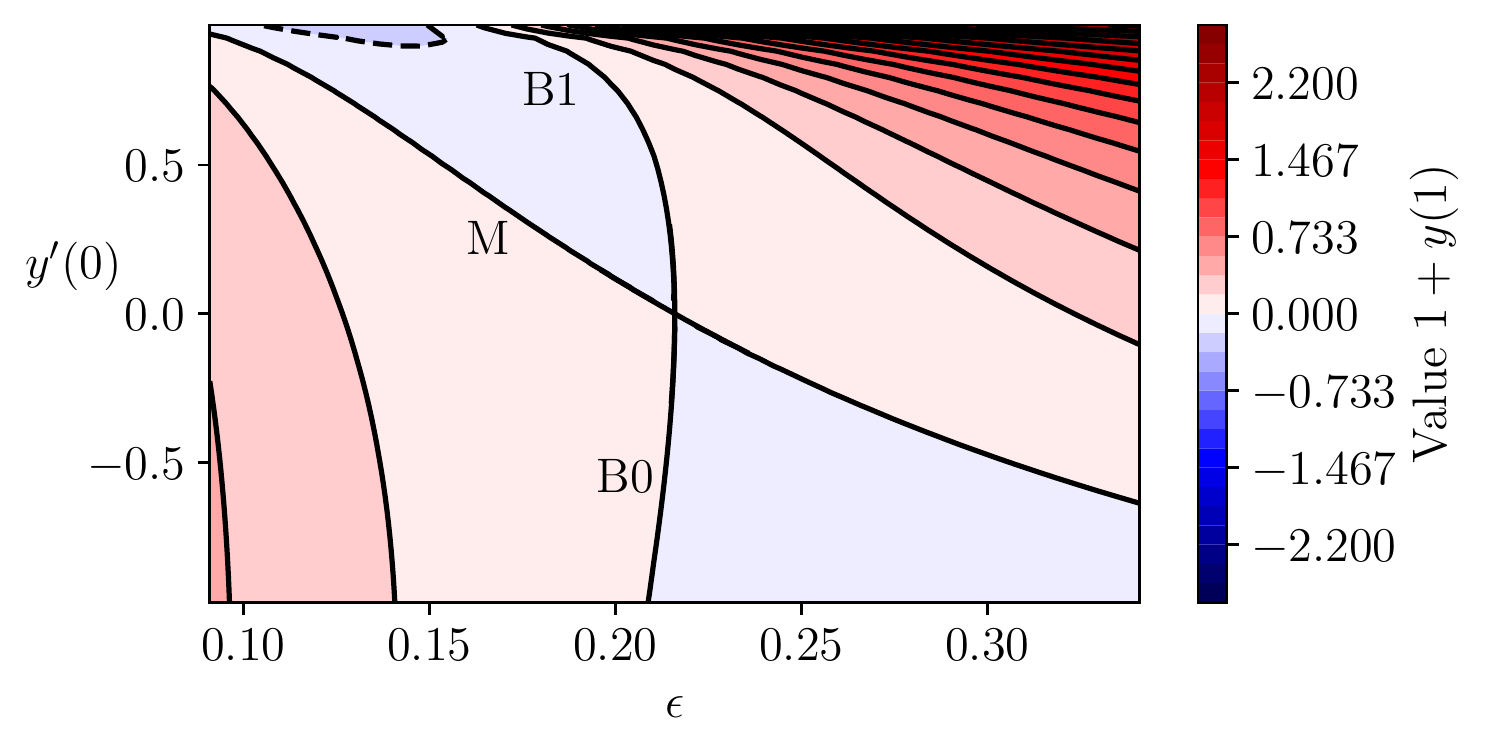}
    \includegraphics[width=4.in]{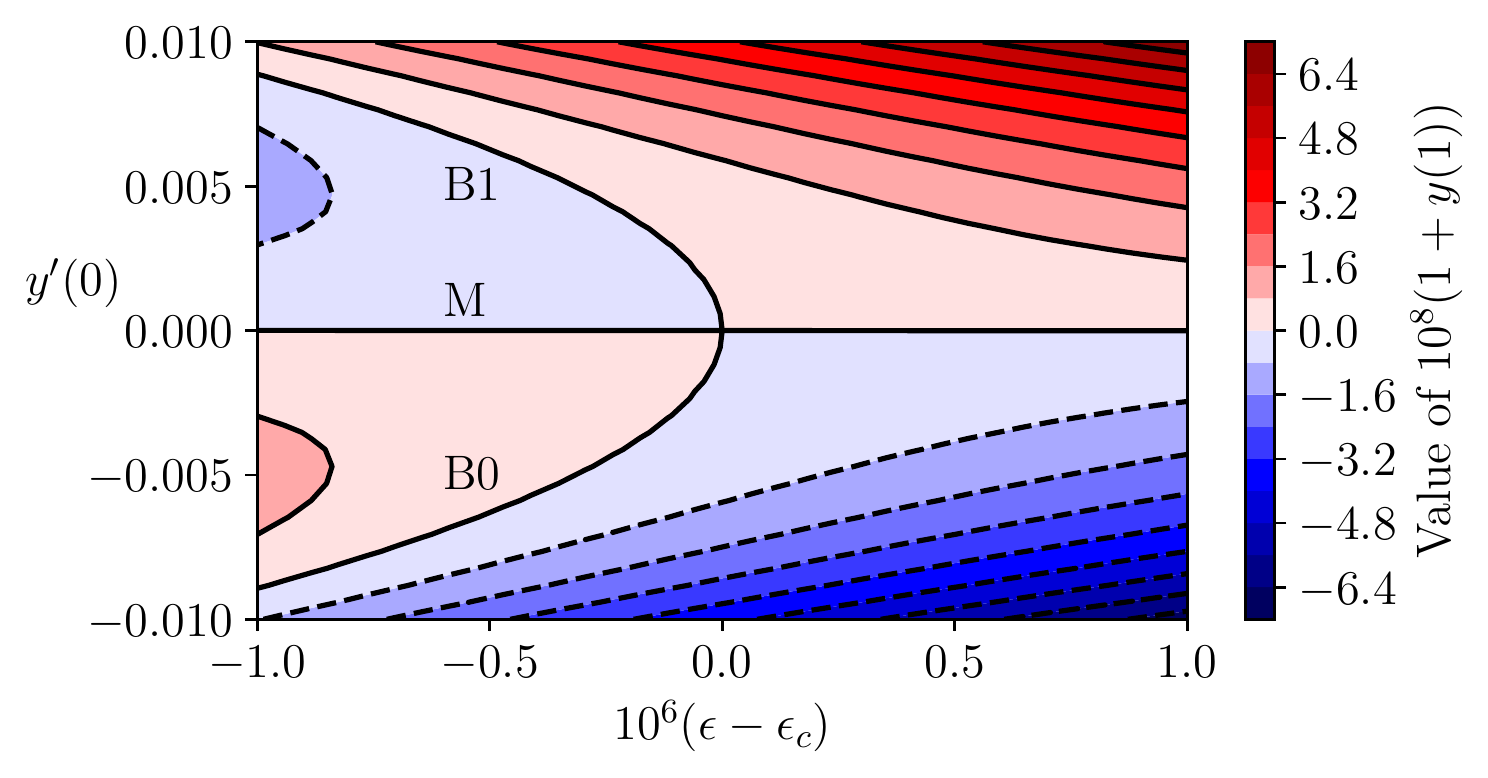}
    \caption{%
      A pitchfork bifurcation occurs in the space of solutions to the boundary-value problem \eqref{eq:ODE_orig}.
      Each point corresponds to integrating the initial-value problem with
      $y(0)=1$, $y'(0)$ determined by the vertical coordinate on the plot, at a
      certain value of $\epsilon$ determined by the horizontal coordinate on
      the plot.  The color denotes the value of $1+y(1)$, with reds denoting
      positive values and blues denoting negative values.  The curves where $1+y(1)=0$ lie between the
      reds and blues; these curves define the set of all $(\epsilon, y'(0))$ that yield solutions of the boundary-value problem. For each $\epsilon < \epsilon_{c}$, there are
      three initial values satisfying $1+y(1)=0$, corresponding to the B0, M,
      and B1 solutions, which confluence at $\epsilon_{c} =
      0.2159869288903\ldots$ and $y'(0)=0$.  The zoom in the lower panel
      expands the neighborhood of the bifurcation point, showing that it has
      the form of a cubic, as expected from the normal form of a pitchfork bifurcation.
    }
    \label{fig:pitchfork}
  \end{figure}%
}

\begin{document}

\maketitle

\begin{abstract}
We revisit a textbook example of a singularly perturbed nonlinear boundary-value problem. Unexpectedly, it shows a wealth of phenomena that seem to have been overlooked previously, including a pitchfork bifurcation in the number of solutions as one varies the small parameter, and transcendentally small terms in the initial conditions that can be calculated by elementary means. Based on our own  classroom experience, we believe this problem could provide an enjoyable workout for students in courses on perturbation methods, applied dynamical systems, or numerical analysis.
\end{abstract}

\begin{keywords}
  boundary layer, singular perturbation, dynamical systems, higher-order matching, matched asymptotics, shooting method
\end{keywords}

\begin{AMS}
  34B16, %
  34E05, %
  34E15, %
  37M20, %
  65L11. %
\end{AMS}

\section{Introduction}
\label{sec:problem}

In many parts of mathematics, physics, engineering, and the life sciences, researchers have developed ingenious techniques for gaining insight into difficult problems  by exploiting the presence of a small parameter in them. These techniques, known as perturbation methods, have shed light on all sorts of fascinating phenomena in fluid dynamics, mathematical biology, optics, chemical engineering, quantum mechanics, plasma physics, climate science, and many other disciplines. See \cite{bender1999advanced, hinch_1991, holmes1995introduction,  kevorkian2012multiple, kuehn2015multiple, lagerstrom2013matched, omalley1991singular} for just a few of the  textbook introductions to perturbation methods and their applications. 

This is the story of a textbook problem that has surprised us over and over again, and for which we have come to feel genuine affection. The problem is to solve the following nonlinear differential equation, subject to the given boundary conditions: 
\begin{align}
  \label{eq:ODE_orig}
  \epsilon y'' &= y y' - y
  \,, &
  y(0) &= 1
  \,, &
  y(1) &= -1
  \,.
\end{align}
Here we use the notation $y' = dy/dx$, and we want to solve for $y(x)$ on the
domain $0 \le x \le 1$ under the assumption that the parameter $\epsilon$ is small: $0 < \epsilon \ll 1$.

We first met problem~\eqref{eq:ODE_orig} in the classic textbook by
Mark Holmes~\cite{holmes1995introduction}. Soon after that book appeared in 1995, one of us (Strogatz) decided to
discuss~\eqref{eq:ODE_orig} in an introductory course on asymptotics and perturbation methods. It quickly became
clear the problem contains unexpected subtleties and riches.

In the years since then, whenever Strogatz has had a chance to teach that
course, he has revisited problem~\eqref{eq:ODE_orig} and always learned something new about it, thanks to the questions and insights of the students, postdocs, and colleagues in attendance, most notably the co-authors on this paper. Together we think we may have finally gotten to the bottom of it.

But because this is meant to be an education paper, we will also be pointing out some of the false turns we took along the way. Confusion is a natural part of doing mathematics. We did not land on the right way to think about~\eqref{eq:ODE_orig} initially; it required lots of trial and error and a willingness to be open and vulnerable about what puzzled us at any given time. Experienced researchers know this is all part of the  process, but we mention it for the sake of students who may be misled by the confident presentations and pristine appearance of the mathematics they see in most textbooks and journal articles. We want to be a little more honest here about how the sausage is actually made.

Our analyses have relied on three parts of applied mathematics:
perturbation theory, nonlinear dynamics, and numerical analysis. As
such, we think problem~\eqref{eq:ODE_orig} could be useful to students
or teachers in any of those subjects. We assume that the reader is
comfortable with asymptotics and perturbation methods at the level of
the books by Holmes~\cite{holmes1995introduction} or Bender and
Orszag~\cite{bender1999advanced}, as well as nonlinear dynamics at the level of the book by Strogatz~\cite{StrogatzNonlinear}. Not much
exposure to numerical analysis is required; a basic knowledge of what
it means to solve an ordinary differential equation  numerically should be sufficient, say at the level of someone who knows how to use \texttt{Mathematica} or \texttt{Matlab} to solve an initial-value problem. For readers who want to immerse themselves in the details, we also provide supplementary notebooks: a \texttt{jupyter} notebook for the numerical methods, and a \texttt{Mathematica} notebook for the analytical calculations~\cite{CodeSupplements}.

Although we have grown enamored of problem~\eqref{eq:ODE_orig} for its  pedagogical value, it would be even more appealing if it also had some real-world applications. Alas, we have not found any so far.  But it does have some close relatives of scientific interest. For example, if we change the sign of the right-hand side of the differential equation in problem~\eqref{eq:ODE_orig}, we get an equation known as the Lagerstrom-Cole equation~\cite{kevorkian2012multiple, lagerstrom2013matched}, which has been studied in aerodynamics and fluid dynamics as a model problem in connection with shock layers~\cite{lagerstrom2013matched}. The differential equation in \eqref{eq:ODE_orig} also pops up in the study of delicate nonlinear phenomena known as canards \cite{benoit1981chasse,  diener1984canard, dumortier1996canard, eckhaus1983relaxation, jardon2020fast, krupa2001extending, kuehn2015multiple}. As defined by Krupa and Szmolyan~\cite{krupa2001extending}, “a canard solution is a solution of a singularly perturbed system which is contained in the intersection of an attracting slow manifold and a repelling slow manifold.” In applications, canards often arise in the analysis of nonlinear oscillations having both fast and slow time scales, as seen in chemistry, neurobiology, electronics, and many other fields~\cite{kuehn2015multiple}.

\section{The surprises}

\figExampleSols

So what are some of the surprises in problem~\eqref{eq:ODE_orig}? The first is how many solutions it has. Holmes~\cite{holmes1995introduction} argued that it has a unique solution, but it turns out that it actually has three. Figure~\ref{fig:example_sols_and_phase_plane} shows their graphs. We call these solutions B0, M, and B1, with the names chosen to indicate that they have a boundary layer at $x=0$, an interior layer in the middle, or a boundary layer at $x=1$. 

The graph of the B0 solution has a very negative initial slope $y'(0)$
at the left endpoint of the domain, while the other two solutions have
$y'(0)$ \emph{extremely} close to 1. In fact, those slopes differ from
1 by a ``transcendentally small'' term, meaning a term that goes to zero
faster than any positive power of $\epsilon$ as $\epsilon \rightarrow
0^{+}$. A common example is an exponentially small term of the form $\exp(-c/\epsilon)$ where $c>0$ is a constant. Normally, such transcendentally small terms are beyond the
reach of basic perturbation theory (although sometimes they can be
handled by more sophisticated techniques known as exponential asymptotics, superasymptotics, hyperasymptotics,
 or asymptotics beyond all
orders~\cite{boyd1998weakly, boyd1999devil}). What especially surprised  
us about problem~\eqref{eq:ODE_orig} is that the leading-order
asymptotics of the initial slope $y'(0)$ could be found by elementary means, even
when the difference between $y'(0)$ and 1 is transcendentally small. 
As we show in section \ref{sec:mapping-from-bvp}, $$1- y'(0) \sim \frac{9}{2\epsilon}\exp\left(-\frac{5}{8\epsilon}\right)$$ for the M solution, while for B1 it is given by $$1- y'(0) \sim \frac{24}{\epsilon}\exp \left(-\frac{3}{2\epsilon}\right).$$
These transcendentally small terms can be calculated using nothing more than higher-order matching, phase plane analysis, and a constant of motion for the associated flow. Not only were we surprised by this good fortune; we were also delighted and relieved by it, given that none of us knew anything about the more sophisticated techniques mentioned above!

Yet another surprise is the behavior of the solutions with respect to $\epsilon$. In section~\ref{sec:pitchfork} we'll see that as we increase $\epsilon$ from zero, the three solutions persist until a critical value $\epsilon_c \approx 0.216$, at which point they collide in a pitchfork bifurcation. For larger values of $\epsilon$, problem~\eqref{eq:ODE_orig} has a unique solution that resembles the M solution above.

\section{Locating the layers: Holmes's analysis}

Our initial thinking about problem~\eqref{eq:ODE_orig} was greatly influenced by what Holmes~\cite{holmes1995introduction} had to say about it. So we summarize his analysis here. We are not going to show the mathematical details yet, because they will appear in later sections. For now we just want to emphasize a question that arises whenever one confronts a differential equation like \eqref{eq:ODE_orig} with a small parameter $\epsilon$ multiplying its highest derivative: 

\vspace{2mm}

\textsc{Question 1:} Does problem~\eqref{eq:ODE_orig} have
any layers (meaning regions of rapid variation in $y$ or its
derivatives), and if so, where are they located?

\vspace{2mm}

This issue comes up midway through Holmes's chapter on matched asymptotic expansions (chapter 2 in Ref.~\cite{holmes1995introduction}). Earlier in the chapter he has already introduced the basic ideas of inner and outer solutions, boundary layers, matching, and uniformly valid composite solutions. There is also a discussion of higher-order matching and examples with multiple boundary layers. In all these prior examples, the layers occur at the endpoints of the domain; in other words, they are genuine ``boundary layers.'' In contrast, problem~\eqref{eq:ODE_orig} is offered as the first
instance of a problem having an ``interior layer.'' 
Holmes writes:

\vspace{1mm}

\begin{quote}
  Generally, when one first begins trying to solve a problem it is not
  known where the layer(s) is. If we began this problem as we did the
  previous two and assumed there is a boundary layer at either one of
  the endpoints, we would find that the expansions do not match. This
  is a lot of effort for no results, but fortunately there is a
  simpler way to come to the same conclusion.
\end{quote}

\vspace{1mm}

\noindent He then offers a convexity argument to rule out boundary layers at either $x=0$ or $x=1$, but is careful to note that ``these are only plausibility arguments and they do not prove anything. What they do is guide the analysis and hopefully reduce the work necessary to obtain the solution.'' 

\subsection{Holmes's plausibility argument}

Holmes first considers a possible boundary layer at $x=0$. He looks at the governing equation $\epsilon y'' = y y' - y$ along with its left boundary condition, $y(0)=1$. Then he sketches the graph of a candidate solution that looks  like our B0 solution in Figure~1, except that he also assumes that $y''> 0$ in the boundary layer; in other words, $y(x)$ is assumed to be concave up for all small $x>0$. On the other hand, although $y''$ (hypothetically) has one sign in the layer, $y$ itself clearly changes sign from positive to negative as it drops from $y(0)=1$ at the boundary to a value $y \approx -2$ where it matches the outer solution. With those sign considerations in mind, Holmes~\cite{holmes1995introduction} continues:

\vspace{1mm}

\begin{quote}
    If the solution behaves like the other example problems we have examined, then near $x=0$ it would be expected that $y' < 0$ and $y'' >0$. This means that $\epsilon y'' >0$ but $y (y' - 1)$ is both positive and negative in the boundary layer. This is impossible \ldots. It is possible to rule out a boundary layer at $x = 1$ in the same way.
\end{quote}

\vspace{1mm}

Having ruled out (convex or concave) boundary layers at either end, Holmes~\cite{holmes1995introduction} next considers the
possibility of an interior layer centered at a point $0< x_0 <1.$ He
works out three possible forms for the inner solution (depending on whether a certain constant of integration is positive, negative, or zero) and shows that only one of these can be matched to the outer solution. Then he uses a symmetry argument to conclude that the layer must be centered at the midpoint of the domain, $x_0=1/2$, and demonstrates that the inner solution must have odd symmetry about that point. Finally he matches the inner and outer solutions, constructs a uniformly valid composite solution, and shows it agrees with a numerical solution obtained for $\epsilon=0.01$. In this way, Holmes convincingly demonstrates the existence and properties of the solution we called M above, a solution of \eqref{eq:ODE_orig} with an interior layer in the middle of the domain. 

But what about the possibility of non-convex or non-concave boundary-layer solutions? Recall that the plausibility argument only rules out solutions with $y''$ having strictly one sign in the layer. As Holmes's careful wording suggests, that loophole could potentially allow for sneakier solutions where $y''$ changes sign within a boundary layer.  As we will see in the next section, that little finesse is precisely what allows the solutions B0 and B1 to exist. 

\section{Another approach to locating the layers: Phase plane analysis}

As we worked through Holmes's analytical approach to problem \eqref{eq:ODE_orig}, we began to wonder if it might be helpful to supplement it with a more geometric style of reasoning known as phase plane analysis. After all, the equation appearing in \eqref{eq:ODE_orig},  $\epsilon y'' = y y' - y$, is a nonlinear, second-order, autonomous differential equation, and phase plane analysis is a powerful tool for illuminating how the solutions to such equations behave. Plus, we have to admit, we have more experience with nonlinear dynamics than perturbation theory, so it felt like a more secure way to approach an  unfamiliar problem.

\subsection{Phase portraits}

To recast the problem into the language of dynamics, we replace the independent variable $x$ with $t$, and we think of it as time. Then the dependent variable $y$ becomes a function of time $t$. The advantage of this approach is that it allows us to use our physical intuition about time and motion and the difference between fast versus slow. Abstract solutions to the differential equation turn into easily pictured trajectories of (imaginary) particles moving around in a two-dimensional space known as the phase plane.   

To construct the phase plane, we convert the second-order equation $\epsilon y'' = y y' - y$ into a pair of first-order equations, and then view those as defining a vector field on the plane. We perform the first step by introducing a new dependent variable $z$, defined as $z=y'$, and then we rewrite $\epsilon y'' = y y' - y$ in terms of $z$ to get $\epsilon z' = y z - y$. By solving for $y'$ and $z'$ and placing them on the left hand side of a pair of first-order differential equations, we obtain the following vector field:
\begin{align}
  \label{eq:ODE_reduced}
  \begin{split}
    y' &= z \,, \\
    z' &= \frac{1}{\epsilon} y(z-1) \,.
  \end{split}
\end{align}

Next we interpret \eqref{eq:ODE_reduced} as a dynamical system. From this perspective, the vector $(y',z')$ then tells us the instantaneous ``velocity'' that an imaginary particle at $(y,z)$ would have at time $t$. As the imaginary particle moves around in the $(y,z)$ phase plane, it traces out a trajectory $(y(t),z(t))$, which is the geometric counterpart of a solution $(y(x), y'(x))$ to the original problem  \eqref{eq:ODE_orig}. 

This construction allows us to visualize how the solutions to \eqref{eq:ODE_reduced} behave by imagining how particles move around in the phase plane. There is no need to be quantitatively precise just yet; a  qualitatively correct picture is enough at this stage. By looking at the signs of $y'$ and $z'$ in \eqref{eq:ODE_reduced} and sketching a few vectors in various parts of the $(y,z)$ plane, we are led to the picture shown in the left panel of Figure~\ref{fig:phase_plane_simple}. This picture is called the phase portrait for the system. It shows that there are three qualitatively different types of trajectories for \eqref{eq:ODE_reduced}: parabolic-looking trajectories that flow from left to right above the horizontal line $z = 1$; a straight trajectory that flows from left to right along the invariant line $z = 1$; and periodic trajectories below $z = 1$ that form closed loops, and on which a particle would circulate round and round,  always moving clockwise. The right panel of Figure~\ref{fig:phase_plane_simple} shows the same information  quantitatively. 

\figPhasePlaneSimple

Incidentally, we can prove that the periodic-looking trajectories in Figure~\ref{fig:phase_plane_simple} are truly periodic and are not merely slowly-winding spirals in disguise. There are two standard ways of proving this, so we will not dwell on the details; see section 6.5 and 6.6 in Ref.~\cite{StrogatzNonlinear} for an introduction. Briefly, one way is to  note that \eqref{eq:ODE_reduced} is a ``conservative'' system. To see this, rewrite it as $y'/z'=dy/dz = \epsilon z/[y(z-1)]$ and then separate variables and integrate to obtain $2 \epsilon [z+\log(1-z)] - y^2 = $ constant. This implicit equation can be shown to define closed curves for all $z<1$. Another way is to observe that that \eqref{eq:ODE_reduced} is also a ``reversible'' system: the vector field~\eqref{eq:ODE_reduced} is unaltered  by the change of variables $(x,y,z) \rightarrow (-x,-y,z)$, which corresponds to a time reversal combined with a mirror reflection across the $z$-axis in Figure~\ref{fig:phase_plane_simple}. From this symmetry we can conclude the trajectories lying below $z=1$ are composed of two left/right mirror-image halves that together form a bilaterally symmetric loop. 

\figPhasePlaneThreeSols

\subsection{An important symmetry}

A stronger reversibility symmetry of \eqref{eq:ODE_reduced} is worth noting: Both the  differential equation $\epsilon y'' = y y' - y$ and its boundary conditions $y(0)=1, y(1)=-1$
are left unchanged by the transformation $$(x,y) \to (1-x,-y).$$ Hence the possible solutions of the original boundary-value problem \eqref{eq:ODE_orig} either come in symmetrical pairs $y(x)$ and $\tilde{y}(x)$, where $\tilde{y}(x) = - y(1-x)$, or the pair degenerates to a single solution with this symmetry, $y(x)=-y(1-x)$. We already saw a visual manifestation of this symmetry in Figure~\ref{fig:example_sols_and_phase_plane}, where B0 and B1 form a symmetric pair and M is self-symmetric. Likewise, when those same solutions are plotted in the phase plane shown in Figure~\ref{fig:phase_plane_three_sols}, the red and blue curves are
paired under the symmetry, and the green curve is self-symmetric. 

\subsection{Boundary conditions}

What about the boundary conditions $y(0)=1$ and $y(1)=-1$ in \eqref{eq:ODE_orig}? How do they enter the phase plane picture? Well, the condition $y(0)=1$ means that at time $t=0$ our imaginary particle must start somewhere on the vertical line $y=1$ in the phase plane (Figure~\ref{fig:phase_plane_three_sols}). Its $z$ coordinate on that line, however, is unspecified and remains to be determined; indeed, the key to solving \eqref{eq:ODE_orig} is to figure out the initial value of $z$ that will enable the moving particle to satisfy the \emph{other} boundary condition, $y(1)=-1$. In dynamical terms, this other boundary condition $y(1)=-1$ is a final condition, not an initial condition. It says that the particle must reach the vertical line $y=-1$ in the $(y,z)$ plane after exactly one unit of travel time.  

Thus, we see what a difficult challenge our imaginary particle is facing. It must find exactly the right place to start on the line $y=1$, such that after it gets carried along by the flow determined by the vector field, somewhat like a tiny speck of leaf being carried downstream by a gentle brook, it manages to land somewhere on the line $y = -1$ precisely when the clock strikes time $t=1.$ 

Figure~\ref{fig:phase_plane_simple} immediately implies that the trajectories on the line $z=1$ or above it are disqualified as candidate solutions because there's no way they can satisfy the boundary conditions: A particle starting on any of them would move monotonically to the right, and so would have no chance of making the leftward journey required to get from $y = 1$ to $y = -1$. That leaves the closed loops below $z=1$ as our only hope. And indeed, we can imagine a particle starting somewhere on the line $y=1$ and then flowing along an arc on one of the closed loops such that if everything is chosen just right (i.e., we pick the right loop to start on), the particle will reach $y = -1$ after exactly one unit of travel time. Figure~\ref{fig:phase_plane_three_sols} shows three arcs that do the trick. 

\subsection{Slow-fast structure}

So far we have not used the assumption that $\epsilon$ is small, but now we will. For $\epsilon \ll 1$, we see from \eqref{eq:ODE_reduced} that the vector field has large regions where the flow is very fast in the vertical direction, with vertical velocities $z'$ of order $O(1/\epsilon)$ occurring at all points in the $(y,z)$ plane  where $z-1 = O(1)$. This region of fast variation, as we will soon see, corresponds to the ``inner region'' in a perturbation treatment via boundary-layer theory. In the phase portrait, it consists of all points outside the thin gray strip shown schematically in Figure~\ref{fig:phase_plane_simple}. The strip does not have well-defined edges, but its blurriness does not matter; the key thing is that its thickness is $O(\epsilon)$, however we define it. (To check its thickness, observe that if  $y$ is $O(1)$ and $z-1 = O(\epsilon)$, then $y'$ and $z'$ are both $O(1)$ in \eqref{eq:ODE_reduced}, indicating the flow is slow compared to the $O(1/\epsilon)$ speeds achieved everywhere outside the strip.) This strip in the phase plane where the motion is comparatively slow corresponds to the ``outer region'' in a boundary-layer treatment. 

It was by contemplating the slow-fast structure of the flow that we originally came to suspect that there might be more than one solution to problem~\eqref{eq:ODE_orig}. As we reconsidered the solution discussed by Holmes, with its interior layer centered at $x = 1/2$, we pictured it as a particle moving with a slow-fast-slow trajectory in the phase phase, as shown in the middle panel of Figure ~\ref{fig:phase_plane_three_sols}. Our imaginary particle spends about half of its travel time dawdling through the initial slow region at the top of the green arc in Figure~\ref{fig:phase_plane_three_sols}, then rockets down and around and up again through the fast region in almost no time at all, and then dawdles through the remaining slow region for the remaining half of its travel time. Why, we wondered, couldn't a particle spend nearly \emph{all} its time in a slow region at the beginning? Or at the end? If the particle started sufficiently close to the invariant line $z=1$, or ended up near there, it seemed like these sorts of solutions should also be possible. 

This intuition turned out to be correct: such solutions do exist. The blue and red curves in Figure~\ref{fig:phase_plane_three_sols} show what they look like as trajectories.

\subsection{More backstory: Puzzling over the initial slope}

It took us considerable trial and error to find these boundary-layer solutions numerically the first time we looked for them in the computer, more than a decade ago. They eluded us completely when $\epsilon$ was very small. Fortunately, for $\epsilon$ only moderately small it was not difficult to find them. For $\epsilon = 0.1$, for example, we found that $y'(0) \approx 0.9999$ yielded a trajectory that was slow at the beginning and fast at the end (the B1 solution), whereas $y'(0) \approx -10.6942$ gave a trajectory that was fast and then slow (the B0 solution). 

The strikingly small difference between 0.9999 and 1 made us wonder what the formula for $y'(0)$ as a function of $\epsilon$ might be for the B1 solution. Likewise, given the size of $-10.6942$, we were curious how negative $y'(0)$ might get for smaller values of $\epsilon$ as we continued tracking the B0 solution. Our numerics couldn't answer these questions at the time, since some of us were naive about computational methods, so we did not know how to solve the boundary-value problem reliably for $\epsilon \ll 1$.

Perhaps the initial slope could be found by asymptotic analysis? We felt sure it could but did not immediately see how to do it. We pose that as our next big question and call it the \emph{puzzle of the initial slope}:

\vspace{2mm}

\textsc{Question 2:} For the three solutions of  problem~\eqref{eq:ODE_orig}, how do their initial slopes $y'(0)$ depend on $\epsilon$ for $0 < \epsilon \ll 1$?

\vspace{2mm}

We will work toward answering that question in the next few sections. But before we leave the phase plane, we should notice that it tells us one more thing of interest. It reveals that the B0 and B1 solutions sneak through the loophole in Holmes's plausibility argument by having layers that are \emph{non-convex}. To see how that conclusion follows from the phase portrait, observe that on any of the trajectories shown in Figure~\ref{fig:phase_plane_three_sols}, the value of $z'$ (the vertical velocity) evidently \emph{changes sign} as the particle goes down and then back up on its journey through the fast region. Since $z'=y''$, that change of sign means the concavity of $y(x)$ changes sign in the layer!

\section{Perturbation theory}
\label{sec:bound-layer}

In this section we solve problem~\eqref{eq:ODE_orig} for $\epsilon \ll 1$, both inside and outside the
boundary layers or interior layers. Then we match the inner and outer solutions and find composite solutions that give uniformly valid asymptotic approximations of $y(x)$ over the whole
domain $0 \le x \le 1$. We perform the match to first order in $\epsilon$ (i.e., we go beyond the leading order of perturbation theory) because it turns out we need this higher-order information to solve the puzzle of the initial slope (Question 2). In that sense, the following analysis provides a motivational case study of why one would ever want to do higher-order matching. 
The details of this analysis are included in the supplementary
\texttt{Mathematica} notebook~\cite{CodeSupplements}.

\subsection{Outer solution}

First we consider the outer region, where regular perturbation
theory applies.  In this region, we expand $y(x,\epsilon)$ in the regular perturbation series
\begin{align}
  y(x, \epsilon) \sim y_{0}(x) + \epsilon y_{1}(x) + O(\epsilon^{2})
  \,,
\end{align}
insert this into the differential equation $\epsilon y'' = y y' - y$ appearing in \eqref{eq:ODE_orig}, and collect terms having like powers of $\epsilon$. At leading order we find
\begin{align}
  y_{0}\, y_{0}' - y_{0} = 0
  \,.
\end{align}
This equation has two possible solutions: $y_{0} = 0$ (which cannot
satisfy the boundary conditions), or $y_{0}' = 1$, yielding $y_{0}(x) = x + a$ for
some real constant $a$. In fact, a bit of study shows that the higher
corrections all satisfy $y_{n}'=0$ for $n>0$, and thus $y_{n} = a_{n}$
for some constants $a_{n}$.  But the leading-order solution $y_{0}(x)$  already fixes the constant by satisfying the boundary condition, so all the higher constants vanish: $a_{n}=0$ for all $n>0$. Therefore $y_0(x)=x+a$ is not merely the zeroth-order approximation to the outer solution; it is the outer solution at \emph{all orders} of $\epsilon$. We can also reach this conclusion by noting that $y_0(x) = x + a$ satisfies the original differential equation $\epsilon y'' = y y' - y$ exactly for all values of $\epsilon$. 

For an outer solution that includes $x=0$ in its domain, we can determine the constant $a$ by applying the boundary condition at that endpoint, and similarly for an outer solution that includes $x=1$.  These
two potential outer solutions that satisfy either the left or right boundary condition are
\begin{align}
  y_{0}^{L}(x) &= 1+x
  \,, &
  y_{0}^{R}(x) &= x-2
  \,.
\end{align}

\subsection{Inner equation}

Now we move on to the inner solutions.  Suppose there is a
layer at $x=x_{0}$. Holmes~\cite{holmes1995introduction} shows that $x_0=1/2$ is the only possible location for an interior layer, and we are going to show that boundary layers can occur at $x_{0}=0$ or $x_0=1$ as well. As before, we refer to these three cases as M (layer in the middle), B0 (layer on the boundary at $x=0$), and B1 (layer on the boundary at $x=1$).

Introduce a layer thickness $\delta$, which is a function of
$\epsilon$ to be determined.  Let
\begin{align}
  X \equiv \frac{x-x_{0}}{\delta} \,,
\end{align}
be a scaled independent variable that describes positions in the layer, and let
\begin{align}
  Y(X) \equiv y(x) = y(x_{0} + \delta X)
  \,
\end{align}
be a new dependent variable that describes how $y$ varies in the layer. Derivatives of the new variable are
\begin{align}
  \frac{dY}{dX} & = \delta y'
  \,, &
  \frac{d^2 Y}{dX^2} &= \delta^{2} y''
  \,.
\end{align}
Now our original differential equation $\epsilon y'' = y y' - y$  becomes
\begin{align}
  \frac{\epsilon}{\delta^{2}}\frac{d^2 Y}{dX^2} = Y \left( \frac{1}{\delta} \frac{dY}{dX} - 1 \right)
  \,.
\end{align}
If $\delta$ is chosen correctly, then $Y$ and all its derivatives should be $O(1)$ as $\epsilon \rightarrow 0$. Thus we find a distinguished limit when $\epsilon/\delta^{2} = 1/\delta$,
or simply $\delta=\epsilon$.  Therefore the inner equation in the layer is given by 
\begin{align}
  \label{eq:ODE_boundary}
  \frac{d^2 Y}{dX^2} = Y \, \frac{dY}{dX} - \epsilon Y
  \,.
\end{align}

\subsection{Leading-order inner solution}
To solve the inner equation \eqref{eq:ODE_boundary} asymptotically, we expand $Y$ as 
\begin{align}
\label{eq:innerexpansion}
  Y(X,\epsilon) \sim Y_{0}(X) + \epsilon Y_{1}(X) + O(\epsilon^{2})
  \,.
\end{align}
Inserting this series into \eqref{eq:ODE_boundary} yields, at leading order,
\begin{align}
 \frac{d^2 Y_{0}}{dX^2}  = Y_{0} \, \frac{dY_{0}}{dX} = \frac{d}{dX} \left( \frac{1}{2} Y_{0}^{2} \right)
  \,.
\end{align}
This can be integrated to obtain 
\begin{align}
  \label{eq:Y0_sol}
  \frac{dY_{0}}{dX} = \frac{1}{2} Y_{0}^{2} + A
  \,,
\end{align}
where $A$ is an integration constant. 

This result has a nice geometrical interpretation. If we recall that $y'=z$ in the phase plane, then \eqref{eq:Y0_sol} shows that, to leading order, the trajectories in the $(y,z)$ plane follow \emph{parabolic} arcs as they move through the inner region where the motion is fast. From the phase plane pictures shown earlier, we know that the only parabolic arcs of interest are those with a negative $z$-intercept, as these are the arcs that lie on the closed loops. Hence
we see that $A$ should be negative, say $A=-\frac{1}{2}b^{2}$.

Then separating the variables in \eqref{eq:Y0_sol} and integrating gives
\begin{align}
  \frac{dY_{0}}{Y_{0}^{2}-b^{2}} &= \frac{1}{2} dX \,, \\
  Y_{0}(X) &= b \tanh \left( c - \frac{b}{2} X \right)
  \,,
  \label{eq:Y0_innersol}
\end{align}
with an additional integration constant $c$.  The two constants $b,c$
are determined by matching to the solution in the outer region.

\subsection{Zeroth-order matching for the symmetric solution M}

Let's see how the matching goes for the interior layer at
$x_{0}=1/2$.  We know we have a symmetric solution that
satisfies $y^{M}(1/2) = 0 = Y^{M}_{0}(0)$ (we could also learn this from the
matching alone).  This symmetry condition tells us that $Y^{M}_{0}$ is an odd function of $X$ and hence $c=0$.  To get the value of $b$,
we need to look at the large-$X$ asymptotic behavior of the inner solution and match it to the asymptotic behavior of the outer solution as $x$ approaches $x_{0}=1/2$ from either side. Thus we need to take the following limit,  $\lim_{X\to \pm \infty} Y^{M}(X)$, and match it to
$y_{0}^{L}(1/2) = 3/2$ and
$y_{0}^{R}(1/2) = -3/2$.  If we recall that $\tanh z$ is odd and 
$\lim_{z\to+\infty} \tanh z = 1$, we see that both limits agree if and only if $b= \pm 3/2$. For either choice of $b$,  our inner solution at zeroth order becomes 
\begin{align}
  \label{eq:Y0M_sol}
  Y^{M}_{0} = - \frac{3}{2} \tanh \left( \frac{3}{4} X \right)
  \,.
\end{align}

Finally we can construct a composite solution $y_c$ by the usual recipe:  $y_{c} = y_{\textrm{outer}} +Y_{\textrm{inner}} - y_{\textrm{match}}$. Carrying out those steps for the zeroth-order approximation to M yields 
\begin{align}
  \label{eq:yMc0}
  y^{M}_{c,0} = x - \frac{1}{2} - \frac{3}{2} \tanh
  \left[
    \frac{3}{4\epsilon} (x-\tfrac{1}{2})
  \right]
  + O(\epsilon)
  \,.
\end{align}

\subsection{Zeroth-order matching for the asymmetric solutions B0 and B1}

We can proceed similarly for the B0 and B1 cases.  In fact we only
need to do the work for one of them, since we can get the other one
from the symmetry transformation $(x,y) \to (1-x, -y)$.  Let's focus on the B0
case, which has a layer at $x_{0}=0$.  Its zeroth-order inner solution \eqref{eq:Y0_innersol}  needs to be matched to
$y_{0}^{R}(0) = -2$, which tells us that $b=-2$.  We also need to
satisfy the boundary condition $y(0)=1$, which is now inside the layer, so $Y_{0}^{B0}(0)=1$
determines the value of $c$ as
$c = -\tanh^{-1}\tfrac{1}{2}$. Hence 
\begin{align}
  Y^{B0}_{0} = -2 \tanh \left( X - \tanh^{-1}\tfrac{1}{2} \right)
  \,.
\end{align}
Similarly for B1 we get
\begin{align}
  Y^{B1}_{0} = -2 \tanh \left( X + \tanh^{-1}\tfrac{1}{2} \right)
  \,.
\end{align}
Now constructing the leading-order composite solutions, we get
\begin{align}
  \label{eq:yB0c0}
  y^{B0}_{c,0} &= x - 2\tanh\left(\frac{x}{\epsilon} - \tanh^{-1}\tfrac{1}{2}\right) + O(\epsilon)
  \,, \\
  \label{eq:yB1c0}
  y^{B1}_{c,0} &= x - 1 - 2\tanh\left(\frac{x-1}{\epsilon} + \tanh^{-1}\tfrac{1}{2}\right) + O(\epsilon)
  \,.
\end{align}

\subsection{Remarks about the zeroth-order composite solutions}

The leading-order composite solutions for M, B0, and B1 look extremely similar to the numerical solutions plotted in
Figure~\ref{fig:example_sols_and_phase_plane}.  They are also uniformly valid,
as can be checked by looking at the difference between the numerical
and analytical solutions, as plotted in Figure~\ref{fig:uniform_conv} for B0. 
Moreover, the layers for all of the zeroth-order solutions are non-convex, as we expected from our earlier phase plane analysis. For example, the B0 solution \eqref{eq:yB0c0} has an inflection point at $x= \epsilon \tanh^{-1}\tfrac{1}{2}$.

\figUniformError

Yet informative as these leading-order solutions are, they are not accurate enough to allow us to calculate the initial slope $y'(0)$ correctly.  So, let's proceed to the next order.  

\subsection{First-order matching for the symmetric solution M}

As we noted earlier,
the outer solutions are $y_0(x) = x+a$ to all orders in $\epsilon$; only
the inner solution and the matching change as we proceed to higher orders. To study the first-order correction $\epsilon Y_{1}$ in the inner solution, we insert the series \eqref{eq:innerexpansion} into our inner equation \eqref{eq:ODE_boundary} and collect the first-order terms. We find that $Y_{1}(X)$ satisfies
\begin{align}
  \frac{d^2 Y_{1}}{dX^2} - \frac{dY_{1}}{dX} \, Y_{0} - \frac{dY_{0}}{dX} \, Y_{1} + Y_{0} = 0
  \,.
\end{align}
Compared to the asymmetric B0 and B1 solutions, the symmetric solution M yields the simplest expressions for $Y_1$, so we focus on that calculation now and relegate the others to Appendix~\ref{sec:high-order-asympt}.

Using the zeroth-order solution for
$Y_{0}^{M}$ from \eqref{eq:Y0M_sol}, we find that $Y_{1}^{M}$ satisfies the following second-order inhomogeneous linear differential equation with variable coefficients: 
\begin{align}
  \label{eq:Y1M_ODE}
  \frac{d^2 Y^{M}_{1}}{dX^2} + \frac{3}{2} \tanh\left(\frac{3X}{4}\right) \, \frac{d Y^{M}_{1}}{dX}
  - \frac{9}{8} \sech^{2}\left(\frac{3X}{4}\right)  \, Y_{1}^{M} = \frac{3}{2} \tanh\left(\frac{3X}{4}\right)
  \,.
\end{align}
We must solve this beast subject to the side condition that
$Y_{1}^{M}(0)=0$ (a condition that follows from symmetry, as $Y^{M}$ is an odd function).  Impressively, Mathematica obliges and produces a long 
expression that we give in \eqref{eq:Y1M} of
Appendix~\ref{sec:high-order-asympt}.  Of the two free constants of
integration appearing in \eqref{eq:Y1M}, we fix the constant $c_{2}$ by imposing the side condition
$Y_{1}^{M}(0)=0$, giving $c_{2}=\pi^{2}/18$.  To determine $c_{1}$, we
need to perform a first-order match to the outer solution $y_{0}^{R}$ on the right (there's no need to worry about additionally matching to $y_{0}^{L}$, the outer solution on the left; matching on the right automatically takes care of matching on the left, by the odd symmetry noted above). 

To perform the matching on the right, we need to know the
asymptotic behavior of $Y_{1}^{M}$ as $X\to +\infty$.  The relevant 
asymptotics are:
\begin{align}
  &&\sinh z \sim \cosh z &\sim \frac{1}{2} e^{z}
  &&(z\to +\infty)
  \,,& \\
  &&\log(1+z) &\sim z
  &&(z\to 0)
  \,,& \\
  &&\Li_{s}(z) &\sim z
  &&(z\to 0)
  \,,&
\end{align}
where $\Li_{s}(z)$ is a special function called a ``polylogarithm'' of order $s$~\cite{NIST:DLMF}.  The polylogarithms can be defined by their power series or recursively from an integral,
\begin{align}
  \Li_s(z) = \sum_{k=1}^\infty \frac{z^k}{k^s}
  = \int_0^z \Li_{s-1}(t) \frac{dt}{t}
  \,.
\end{align}
For $s<2$, they are elementary functions, for example $\Li_1(z)=-\log(1-z)$. 

With further help from Mathematica we eventually find
\begin{align}
  Y_{1}^{M} \sim X + \frac{2}{3} (c_{1}-1-\log 4)
  \qquad (X\to\infty)
  \,.
\end{align}
Miraculously---and yet not miraculously at all, if one believes in perturbation theory---this function has the exactly the right large-$X$ behavior, $Y_1^M\sim X$, needed to match onto the \emph{outer} solution, if the constant is correct!
Recall that at zeroth order, $\lim_{X\to\infty} Y_{0}^{M}(X) = -3/2$ already
matches the value of $\lim_{x\to {1/2}^{+}} y_{0}^{R}(x)$.
Therefore we want the additive constant above to vanish (or else we
would make an $O(\epsilon)$ error in the matching). Hence
\begin{align}
\label{Mc_1}
  c_{1} = 1+\log 4
  \,.
\end{align}
Finally, by again invoking the recipe $y_{c} = y_{\textrm{outer}}+Y_{\textrm{inner}}-y_{\textrm{match}}$ for forming a composite solution, we obtain the composite solution up to terms of order $\epsilon^2$. Remarkably, the $y_{\text{match}}$ that we need to subtract off here is simply the exact outer solution $y_0^R$. So the whole composite solution boils down to the inner solution $Y_0^M+\epsilon Y_1^M$. Thus, the first-order composite solution for M is:  
\begin{align}
  y_{c,1}^{M} =&
  -\frac{3}{2} \tanh \left(\frac{3 X}{4}\right)
  \nonumber\\
  &{}+
  \frac{\epsilon}{72}
  \sech^2\left(\frac{3 X}{4}\right)
  \Bigg\{4 \pi^2 + 48 \text{Li}_2\left(-e^{-3
      X/2}\right)+9 X (8+3 X)
  \nonumber\\
  &\qquad +48 \sinh \left(\frac{3 X}{2}\right)
  \log \left(2 \cosh \left(\frac{3 X}{4}\right)\right)
  \Bigg\}
  +O(\epsilon^{2})
  \,,
  \label{eq:yMc1}
\end{align}
where $X=(x-\frac{1}{2})/\epsilon$.

\figUniformErrorHigherOrder

We proceed similarly for the B0 and B1 solutions, collecting the results in Appendix~\ref{sec:high-order-asympt}. Figure~\ref{fig:uniform_conv_higher} confirms that the error between the first-order asymptotic solution for B0 and a numerical solution (presumed to be close to exact) truly does shrink in proportion to $\epsilon^2$, as it should for a first-order match. This sort of test provides a reassuring check when doing complicated numerics and asymptotics.

\section{Solving for the initial slope}
\label{sec:mapping-from-bvp}

Now that we have constructed asymptotic approximations to the three solutions M, B0, and B1 of our original problem \eqref{eq:ODE_orig}, we can use those approximations to estimate the initial slope $y'(0)$ in each case. Knowing this initial slope is theoretically interesting since (as we'll see) it depends on $\epsilon$ in an intriguing way. But it's also practically useful information: having a good approximation of the initial slope helps us solve the boundary-value problem \eqref{eq:ODE_orig} numerically.  

\figShootingTrials

One computational approach to solving a boundary-value problem is called the ``shooting method.''  A more thorough discussion can be found in several textbooks on numerical methods, for example~\cite{NumericalRecipes}.
We also provide a supplementary \texttt{jupyter} notebook which
implements the numerics described in this section~\cite{CodeSupplements}.

Figure~\ref{fig:shooting_trials} shows an example of how shooting applies to our problem. We start a trial solution at $y(0)=1$ and launch it with some initial slope $y'(0)$, somewhat like shooting an artillery shell at an intended target. In our case, the target is the point at the other boundary condition: $x=1, y = -1$. Incorrect choices of the initial slopes at $x=0$ will produce solutions that fail to hit 
the boundary condition at $x=1$. The three curves shown on the left in Figure~\ref{fig:shooting_trials}  indicate what happens if we aim too high or too low or just right. 

If we compute where we hit for \emph{many} possible $y'(0)$, and plot the resulting $y(1)$'s versus $y'(0)$, we get the graph of the target function shown in the right panel of Figure~\ref{fig:shooting_trials}. To find a numerical solution of our problem, then, we just have to figure out where the graph crosses the dashed horizontal line $y=-1$.  This is a standard numerical task; it amounts to a root-finding
problem in a small neighborhood of the solution. Notice that for the value  $\epsilon = 0.1$ used to make Figure~\ref{fig:shooting_trials}, the graph of the target function crosses the dashed line in three places. Those are the desired initial slopes of our three solutions. 

To find analytical estimates of these slopes, we will see next why we needed to go to the trouble of doing higher-order matching.  

\subsection{Initial slope for B0}

Let's start by evaluating the value
and slope of the zeroth-order composite solution for B0:
\begin{align}
  \label{eq:IC_B0_c0}
  y_{c,0}^{B0}(0) &= 1
  \,, &
  y_{c,0}^{B0 \, \prime}(0) &= - \frac{3}{2\epsilon}+1
  \,.
\end{align}
The leading behavior here is $y_{B0}' \sim -3/(2\epsilon)$, and
this is indeed correct, as we will soon see.  However, the slope in
\eqref{eq:IC_B0_c0} has an error that is $O(1)$.  This
can be seen graphically in Figure~\ref{fig:uniform_conv}: Within the
layer of width $O(\epsilon)$, there is an
$O(\epsilon)$ error in the value of $y$, leading to an
$O(1)$ error in the slope.  Or, reiterating the point in
another way: The derivative of a function's asymptotic series (at some
order) is not necessarily the asymptotic series of the derivative of the function
(to the same order).  Since we have the composite solution to higher
order, we can easily check what the correction is by using the results for B0 in Appendix~\ref{sec:high-order-asympt}. Evaluation and differentiation of the first-order composite solution at $x=0$ yield
\begin{align}
  \label{eq:IC_B0_c1}
  y_{c,1}^{B0}(0) &= 1
  \,, &
  y_{c,1}^{B0\,\prime}(0) &=  - \frac{3}{2\epsilon} + 1 + \log 16 + O(\epsilon)
  \,.
\end{align}
This is now the full result for $y_{B0}'(0)$ up to errors of order
$O(\epsilon)$. As shown by the top curve of Figure~\ref{fig:initial_slopes}, this asymptotic result agrees nicely  with the value of $y'(0)$ produced by the shooting method.

\subsection{Initial slope: The tempting but wrong way}

For the B0 case, we saw that the $O(1)$ term in the initial slope changed from 1 to $1+\log 16$ when we went to first order. This change in the $O(1)$ term is a rather small effect as $\epsilon \to 0$,
since the $-3/(2\epsilon)$ term dominates in this limit anyway. But to our  surprise (probably because of our inexperience in these matters), we  soon discovered that the B1 and M
cases are much more subtle: differentiating the first-order composite solution does not give the correct initial slope, not even at leading order in $\epsilon$. 

It's an instructive trap to fall into, so let's take the plunge. For simplicity, let's work with M. If we calculate the initial value and initial slope of the zeroth-order composite
solution $y_{c,0}^{M}$ given in \eqref{eq:yMc0}, we find
\begin{align}
  \label{eq:IC_M_c0}
  y_{c,0}^{M}(0) &= -\frac{1}{2} + \frac{3}{2} \tanh \frac{3}{8\epsilon}
  \sim 1 - 3 e^{-\frac{3}{4\epsilon}} + \text{TST}
  \,, \\
  y_{c,0}^{M \,\prime}(0) &= 1 - \frac{9}{8\epsilon} \sech^{2} \frac{3}{8\epsilon}
  \sim 1 - \frac{9}{2\epsilon} e^{-\frac{3}{4\epsilon}} + \text{TST}
  \,,
\end{align}
where TST stands for ``transcendentally small terms,'' i.e.,\
terms that are smaller than any power of $\epsilon$ times the smallest
reported term.

Should we trust these results?
Let's check by going to the next order of perturbation theory. Our first-order composite solution yields 
\begin{align}
  \label{eq:IC_M_c1}
  y_{c,1}^{M}(0) &\sim 1 - e^{-\frac{3}{4\epsilon}} \left[\frac{3}{8\epsilon} + 4 + O(\epsilon) \right] + \text{TST}
  \,, \\
  y_{c,1}^{M\,\prime}(0) &\sim 1 - e^{-\frac{3}{4\epsilon}} \left[ \frac{9}{16\epsilon^{2}} + \frac{9}{2\epsilon} + \frac{\pi^{2}}{3} + O(\epsilon^{2})  \right]
  + \text{TST}
  \,.
\end{align}
This doesn't look good at all! Going to the next order in $y$ has
resulted in a change at a \emph{lower} order to the values of $y(0)$ and $y'(0)$.
Apparently we can't trust this. So for this problem at least, naively differentiating the composite solution does not give us the correct initial slope.

\subsection{Initial slope: The right way}

Instead of the approach above, let's turn to a more global analysis, using our knowledge of the structure of solutions in the phase plane. The key insight is that we can use a conserved quantity (also known as a constant of motion, or a first integral) to transfer trustworthy information from inside a layer to a distant point outside the layer where we want to calculate an initial slope. For example, we can transfer information from the M layer at $x_0 = 1/2$ all the way over to $x=0$; this trick is how we are going to extract the leading (but still minuscule!) transcendentally small term in M's initial slope. The advantage of using a conserved quantity is that it is exact; it allows us to shuttle information around in the phase plane with perfect fidelity. 

\subsubsection{Conserved quantity for the trajectories}
\label{sec:conserved}

To set the stage to perform the desired transfers of information, we need to gather a few facts about the conserved quantity. 

We have already mentioned that our vector field \eqref{eq:ODE_reduced} is conservative. Recall that at each point in the $(y,z)$ phase plane, 
\begin{align}
  \frac{dz}{dy} = \frac{z'(t)}{y'(t)} = \frac{y(z-1)}{\epsilon z}
  \,,
\end{align}
which can be separated and integrated to yield
\begin{align}
  \label{eq:y_z_fund_rel}
  2 \epsilon [z+\log(1-z)] = y^2-C^2
  \,,
\end{align}
where $C$ is a constant that labels the trajectories. Since $C$ is constant, the quantity $$2 \epsilon [z+\log(1-z)] - y^2$$ remains unchanged as $y(t)$ and $z(t)$ flow along a trajectory and hence is a ``conserved quantity.''

For a given value of $C$, we can generate two explicit formulas for the trajectories as curves in the $(y,z)$ plane. Either we can  write $y(z)$ in the right or left half plane by solving for
$y$ and using the relevant branch of the square root.  Or we can solve for
$z(y)$ if we allow ourselves to use the implicitly-defined Lambert $W$
function, which satisfies
\begin{align}
  \label{eq:W_def}
  W(z)\exp(W(z))=z
  \,.
\end{align}
This equation has
multiple solutions, giving the multiple branches $W_{n}(z)$.
For more on the Lambert $W$ function, see~\cite{corless1996lambertw}.
From \eqref{eq:y_z_fund_rel} the explicit solution for $z$  is
\begin{align}
  \label{eq:z_of_y_explicit}
  z = 1+W_{n}
  \left(
    -\exp
    \left(
      -1 + \frac{y^2-C^2}{2\epsilon}
    \right)
  \right)
  \,.
\end{align}
For $0 \le z < 1$, we want the branch $W_{0}$. For $z \le 0$, we want
the branch $W_{-1}$.

The existence of a conserved quantity for Eq.~\eqref{eq:ODE_reduced} suggests that the dynamical system might actually be a Hamiltonian system in appropriate coordinates. Indeed, Eq.~\eqref{eq:ODE_reduced} is Hamiltonian for $z < 1$: Canonical variables are $Q=y$ and $P=\log(1-z),$ in terms of which the vector field corresponding to~\eqref{eq:ODE_reduced} becomes $Q^\prime = 1-e^P, P^\prime = Q/\epsilon$, with a corresponding Hamiltonian $H(P,Q) = P-e^P-Q^2/(2\epsilon)$.

\subsubsection{Using the conserved quantity to transfer information about slopes}
Having set the stage, we're now ready to explain how to transfer slope information reliably with the help of the conserved quantity $2 \epsilon [z+\log(1-z)] - y^2$. Here's the idea: Recall that, by definition, $z=y'$, so $z$ represents a slope $y'$ on a graph of $y$ versus $x$. Suppose we have one point $(y_{1}, z_{1})$ on a solution where we trust
the slope $z_1$. Then, by using the constancy of $C$, we can use this information to get the slope $z_{0}$ at
some other $y_{0}$ on the same solution. For the M and B1 solutions, we will take
$(y_{1},z_{1})$ to be a convenient point inside the layer where the inner solution is trustworthy and we can
compute $z_1=O(\epsilon^{-1})$ accurately.  Then we will
transfer this information over to the initial condition at $x=0$,
where $y_{0}=1$ and we seek $z_0$.

Taking two copies of \eqref{eq:y_z_fund_rel} on the same curve of
constant $C$ and eliminating $C$, we find that our two points are related
by
\begin{align}
  \label{eq:two_points_same_C}
  2 \epsilon [z_{0}+\log(1-z_{0})] - y_{0}^2 = 2 \epsilon [z_{1}+\log(1-z_{1})] - y_{1}^2
  \,.
\end{align}
Let $y_{0}=1$. Solve for the desired initial slope $z_{0}$ by using the appropriate branch of the Lambert $W$ function:
\begin{align}
  \label{eq:z0_from_other_point}
  1 - z_{0} = - W_{0}
  \left(
    - [1-z_{1}] \exp \left( \frac{1-y_{1}^{2}}{2\epsilon} - (1-z_{1}) \right)
  \right)
  \,.
\end{align}
Equation~\eqref{eq:z0_from_other_point} is the key formula for transferring information from one point to another.

Now if we can find trustworthy values of $y_1$ and $z_1$ to plug in, we'll be in business. Let's illustrate the idea with the B1 solution. Remember that B1 is related to the B0 solution by the symmetry transformation $(x,y,z) \to (1-x,-y,z)$. That means that the \emph{final} slope in B1's layer at $x=1$ is the same as the \emph{initial} slope in B0's layer at $x=0$. But we have already calculated B0's initial slope reliably! It's given in \eqref{eq:IC_B0_c1}. (It's reliable because it was calculated inside a layer, namely, the layer at $x=0$ for the B0 solution.) 

Hence the trustworthy point $(y_1,z_1)$ to pick for the B1 solution is 
\begin{align}
  \label{eq:B1_other_point}
  y_{B1}(1) &= -1
  \,, &
  z_{B1}(1) &= -\frac{3}{2\epsilon} + 1 + \log 16 + O(\epsilon)
  \,.
\end{align}
When we plug this point into \eqref{eq:z0_from_other_point}, we soon discover something fascinating about the initial slope $z_{B1}(0)$: it deviates from 1 by a transcendentally small quantity. To see this, however, takes a few more steps. Upon performing the substitution we first obtain
\begin{align}
  \label{eq:z_B1_with_W}
  1 - z_{B1}(0) = - W_{0}
  \left(
    \left[-\frac{3}{2\epsilon} + \log 16 + O(\epsilon)\right] \exp \left( - \frac{3}{2\epsilon} + \log 16 + O(\epsilon) \right)
  \right)
  \,, 
\end{align}
which may look opaque to anyone  unfamiliar with the Lambert $W$ function. Fortunately this expression can be simplified by using the fact that $W_{0}(z)$ is analytic at $z=0$ and its power series is
convergent within a radius $|z| < 1/e$ (see~\cite{corless1996lambertw}
for a thorough treatment).  The first few terms are
\begin{align}
  W_{0}(z) = z - z^{2} + O(z^{3})
  \,,
\end{align}
as can be verified by substituting back into the defining
equation~\eqref{eq:W_def}.  We only need to keep the first term since
the second term is already transcendentally small in $\epsilon$ 
relative to the first term, in light of the form of the argument of $W_{0}$ in \eqref{eq:z_B1_with_W}. 

Thus, we finally arrive at the correct asymptotic
behavior for B1's initial slope: 
\begin{equation}
  1 - y_{B1}'(0) =
  \left[\frac{3}{2\epsilon} - \log 16 + O(\epsilon)\right] \exp \left( - \frac{3}{2\epsilon} + \log 16 + O(\epsilon) \right)
  + \text{TST}
  \,, \\
  \label{eq:z_B1}
\end{equation}
which can be further simplified to 
\begin{align}
  1 - y_{B1}'(0) =&
  \left[\frac{24}{\epsilon} + O(1)\right] \exp \left( - \frac{3}{2\epsilon} \right)
  + \text{TST}
  \,.
  \label{eq:z_B1_simpler}
\end{align}   
To obtain this last equality, we simplified \eqref{eq:z_B1} by replacing $\exp(\log 16)$ with $16$ and  $\exp(O(\epsilon))$ with $1 + O(\epsilon)$. That relative error at $O(\epsilon)$ in \eqref{eq:z_B1} times the leading $\epsilon^{-1}$ term results in our ignorance of the subdominant $O(1)$ term in the prefactor multiplying the controlling exponential in \eqref{eq:z_B1_simpler}. Nevertheless, we still have enough information to nail down the leading-order $\epsilon$-dependence of the initial slope, given by the term $[24/\epsilon] \exp(-3/2\epsilon)$. 

In retrospect, these error considerations clarify \emph{why} we needed go to the bother of approximating the inner solution with higher-order perturbation theory: It's because we needed the error in the argument of the exponential in \eqref{eq:z_B1} to be $O(\epsilon)$, which necessitated an error of $O(\epsilon^{2})$ in the inner solution. 

\figInitialSlopes

For the M solution, we redo the calculation above, except now we use the slope in the layer at $x=1/2$ as our trustworthy value of $z_1$. We get that slope with sufficient precision by simply evaluating the derivative of the first-order composite solution  \eqref{eq:yMc1} at $x=1/2$.  We also recall that $y$ vanishes at $x=1/2$, by the odd symmetry of M about its midpoint. Hence the trustworthy values of $y$ and its slope are
\begin{align}
  \label{eq:M_other_point}
  y_{M}(1/2) &= 0
  \,, &
  y_{M}'(1/2) &= - \frac{9}{8\epsilon} + 1 + \log 4
  + O(\epsilon)
  \,.
\end{align}
Plugging this phase space point into
\eqref{eq:z0_from_other_point}, we find
\begin{align}
  \label{eq:z_M_with_W}
  1 - y_{M}'(0) &= - W_{0}
  \left(
    \left[-\frac{9}{8\epsilon} + \log 4 + O(\epsilon)\right] \exp \left( - \frac{5}{8\epsilon} + \log 4 + O(\epsilon) \right)
  \right)
  \,, \\
 & =
  \left[\frac{9}{8\epsilon} - \log 4 + O(\epsilon)\right] \exp \left( - \frac{5}{8\epsilon} + \log 4 + O(\epsilon) \right)
  + \text{TST}
  \,, \\
  \label{eq:z_M}
   & =
  \left[\frac{9}{2\epsilon} + O(1)\right] \exp \left( - \frac{5}{8\epsilon} \right)
  + \text{TST}
  \,.
\end{align}
Here, as in \eqref{eq:z_B1}, we have enough precision to  nail the leading-order term after expanding $\exp(O(\epsilon))$, but not enough to determine the subleading $O(1)$ correction.

Figure~\ref{fig:initial_slopes} compares these analytical initial slopes against numerical results obtained from the shooting method. In all three cases, the  agreement between asymptotics and numerics is excellent. And for both M and B1, the agreement extends over more than twenty orders of magnitude. Wow! 

\section{Pitchfork bifurcation}
\label{sec:pitchfork}

As we've seen, the nonlinear boundary-value problem \eqref{eq:ODE_orig} has three solutions when
$\epsilon$ is sufficiently small. The final surprise in this problem comes when we examine what happens to the three corresponding initial slopes,
$y'(0)$, as we increase $\epsilon$ away from 0.  We can mull over
the three equations~\eqref{eq:IC_B0_c1}, \eqref{eq:z_B1}, and
\eqref{eq:z_M}, or examine Figure~\ref{fig:initial_slopes} that
summarizes all three.  For sufficiently small $\epsilon$, we have the
ordering $y_{B0}'(0) < y_{M}'(0) < y_{B1}'(0)$.  As $\epsilon$ increases, $y_{B0}'(0)$ increases (becomes less negative) and moves toward zero, while both 
$y_{M}'(0)$ and $y_{B1}'(0)$ decrease 
away from 1 and also move toward zero. 

This behavior suggests that at sufficiently large $\epsilon$, there is a possibility that two or even three solutions might approach each other and merge. It also suggests that the initial slopes of the merging solutions might lie somewhere close to $y'(0) = 0$.  We first discovered experimentally that this actually happens; a three-way merger occurs through a pitchfork bifurcation precisely when $y'(0) = 0$. After the fact, we were able to establish what analytical conditions describe the bifurcation and identify the critical value $\epsilon_{c}$ where it
occurs.

\figPitchfork

\subsection{Visualizing and explaining the pitchfork bifurcation}

We visualize the merger of solutions by plotting contours of $1+y(1)$ in the
$(\epsilon, y'(0))$ plane in Figure~\ref{fig:pitchfork}.  Satisfying the
boundary-value problem means finding the contours where $y(1)+1=0$, which occurs at the
boundary between reds and blues in the figure.  Below the critical
value $\epsilon<\epsilon_{c}$, there are three solutions; these
bifurcate at $\epsilon=\epsilon_{c}$ and $y'(0)=0$. For large values of $\epsilon$, there is only one solution.  In the broad view shown in the top panel, there is a
rather obscure shape for the solution set.  In the bottom panel, we
zoom in to the sixth decimal place in $\epsilon$, and see a classic
pitchfork shape, which we will explain below.

The first condition for the bifurcation comes from understanding the
relationship between the initial and final slopes for the B0 and B1
solutions, which are paired under the symmetry
$(x,y,z) \to (1-x,-y,z)$.  Let's return to
\eqref{eq:two_points_same_C}, the relationship between two points
$(y_{0},z_{0})$ and $(y_{1},z_{1})$ on the same trajectory.
Take these two points to be the endpoints of the trajectory, where $y_{0}=1$ and
$y_{1}=-1$.  Then the $y^2$ terms cancel, leaving a simpler condition
relating the two endpoint slopes:
\begin{align}
  \label{eq:two_endpoints_slopes}
  z_{0}+\log(1-z_{0}) &= z_{1}+\log(1-z_{1})
  \,.
\end{align}
For convenience, let us define a function
\begin{align}
  f(z)=z+\log(1-z)
  \,,
\end{align}
whose domain is $z<1$.  Equation~\eqref{eq:two_endpoints_slopes} says
that $f(z_{0})=f(z_{1})$.  Now there are two possibilities.  For the
self-symmetric solution M, the slopes at the endpoints agree automatically, so $z_{0}=z_{1}$ and  \eqref{eq:two_endpoints_slopes} is vacuously satisfied. This case tells us nothing new. However, for the asymmetric solutions B0 and B1, there are two distinct slopes, $z_{0}\ne z_{1}$, and we can show that they must have opposite signs. 

This result follows from the shape of the graph of $f$. Note first that $f'(z)=1-1/(1-z)$, so $f$ has only one local
extremum, at $z=0$, where $f(0)=0$.  The second derivative
$f''(z) = -(1-z)^{-2}$ is strictly negative everywhere, so $z=0$ is in fact a global maximum.  For $z<0$, $f(z)$ is monotonically increasing, while for $0<z<1$, $f(z)$ is monotonically decreasing.
Therefore, in the case where we have two
distinct slopes $f(z_{0})=f(z_{1})$ but $z_{0} \ne z_{1}$, we can
deduce that they must lie in the intervals $z_{0} < 0$ and $0 < z_{1} < 1$,
or vice versa, and hence have opposite signs, as claimed. 

Graphically, every horizontal slice through the graph of $f(z)$ below its global maximum cuts it in two places: one value $z_0 < 0$, and the other value at $0<z_1<1$.  As we increase our horizontal slice towards the maximum, these two roots both approach and ultimately confluence at $z=0$. 

Finally let's apply this knowledge to the slopes $z_{B0}$ and
$z_{B1}$, when both solutions exist.  The B0 solution has endpoint
slopes $y'(0)=z_{B0}$ and $y'(1)=z_{B1}$, and the B1 solution
exchanges these two.  By the previous analysis, $f(z_{B0})=f(z_{B1})$
with $z_{B0} < 0$ while $0 < z_{B1} < 1$.  We know that as $\epsilon$
increases, these two slopes approach each other.  The only place they
can merge is at $z_{B0}=z_{B1}=0$, which happens at some
$\epsilon=\epsilon_{c}$.

What about the intermediate M solution?  How does its slope compare to the other two? We know that the three slopes are the 
three roots of the target function plotted in the right panel of
Figure~\ref{fig:shooting_trials}.  Is it possible for $z_{M}$ to
confluence with one of $z_{B0}$ or $z_{B1}$ \emph{before} the outer
two roots meet each other?  No: the M solution has the
same slope at each endpoint, so if it confluences with \emph{either} of
the other two roots, \emph{all three} must confluence simultaneously. Generically, this leads us to expect that in a
small neighborhood of the bifurcation point, the error function $y(1)+1$ plotted in Figure~\ref{fig:shooting_trials} should be approximately a cubic of the form
\begin{align}
  y(1)+1 \approx A(y'(0))^{3} + B (\epsilon-\epsilon_{c}) y'(0)
  \,.
\end{align}
We further expect $A>0$ and $B>0$ by examining the shape of the curve in
Figure~\ref{fig:shooting_trials}.  
For $\epsilon>\epsilon_c$, there is only one real solution; for $\epsilon<\epsilon_c$, there are three real solutions; and as $\epsilon\to\epsilon_c$ from below, these three solutions degenerate into a triple root.  All of these expectations are confirmed by the classic pitchfork scenario seen in the lower panel of Figure~\ref{fig:pitchfork}.

\subsection{Calculating the pitchfork bifurcation value}

The preceding analysis told us that the bifurcation happens when \eqref{eq:ODE_orig} admits a solution with $y'(0) = 0$. We can now use that result to determine $\epsilon_{c}$.  

To do so, we first observe that the solution of \eqref{eq:ODE_orig} when $\epsilon = \epsilon_{c}$  must satisfy the boundary condition $y(0)=1$ as well as the bifurcation condition $y'(0)=0$.  These conditions then uniquely determine the corresponding critical trajectory: it starts at $(y(0),z(0)) = (y(0), y'(0)) = (1, 0)$ and also satisfies  $2 \epsilon [z+\log(1-z)] = y^2-C^2$ by \eqref{eq:y_z_fund_rel}. So by plugging in $y(0)=1$ and $z (0)=0$ we see that $C=1$ on the critical trajectory.  

Next, we use $C=1$ to find two conditions on $\epsilon_{c}$. Both conditions relate $\epsilon_{c}$ to $\zc$, defined as the minimum (i.e.~most negative) value of $z$ on the critical trajectory. By solving those two conditions simultaneously, we find $\zc$ and $\epsilon_{c}$, as follows. 

To obtain the first condition on $\epsilon_{c}$, we note that $z'=0$ when $z(t)$ reaches its minimum along the  trajectory, which implies (from \eqref{eq:ODE_reduced}) that $y=0$ there. Like all other points on the critical trajectory, this point $(y,z) = (0, z_c)$  must satisfy \eqref{eq:y_z_fund_rel} with $C=1$. This gives us our first condition on $\epsilon_{c}$: 
\begin{align}
  \label{eq:eps_crit_z_min}
  2\epsilon_{c}\left[\zc+\log(1-\zc)\right] = -1
  \,.
\end{align}

The second condition is that the time required for the critical trajectory to go from $y=1$ to $y=-1$ is $t=1$, just as it is for every solution of our original boundary-value problem. But for the critical trajectory, we can say more. The critical trajectory is self-symmetric, which implies that the time to go halfway is simply $t=1/2$. And at that halfway point, the trajectory is at the point we have just been discussing, $(y,z) = (0, z_c)$. To translate these observations into the condition we seek, we need to find a formula for the travel time. The trick is to write $dt$ in terms of $z$ on the trajectory and then integrate. Recall that $z' = y(z-1)/\epsilon$, from \eqref{eq:ODE_reduced}, so $dt =dz/z' = \epsilon \, dz/[y(z-1)]$. Thus 
$$\frac{1}{2} = \int_{0}^{\frac{1}{2}} dt = \int_{0}^{z_c} \frac{\epsilon_c \,dz}{y(z-1)}$$ 
where we are thinking of $y$ as a function of $z$ on the critical trajectory. That function can be written explicitly by using the fact that $C=1$ on the critical trajectory; solving $2 \epsilon_c [z+\log(1-z)] = y^2-1$ for $y$ and noting that $y \ge 0$ on the first half of the trajectory, we take the positive square root and obtain  $$y=\sqrt{1+2\epsilon_{c}[z+\log(1-z)]}.$$ Hence our travel  time condition becomes
\begin{align}
  \label{eq:eps_crit_time}
  \frac{1}{2} =
  \int_{\zc}^{0} \frac{\epsilon_{c} \, dz}{(1-z)\sqrt{1+2\epsilon_{c}[z+\log(1-z)]}}
  \,.
\end{align}

Now we have two conditions, \eqref{eq:eps_crit_z_min} and \eqref{eq:eps_crit_time}, in terms of $\zc$ and $\epsilon_{c}$.
We use \eqref{eq:eps_crit_z_min} to eliminate $\epsilon_{c}$,
leaving an integral equation for $\zc$ alone to satisfy.
After a bit of algebra, this combined condition is
\begin{align}
  0 = g(\zc) = \zc+\log(1-\zc) +
  \int_{\zc}^{0} 
  \frac{dz}{(1-z)\sqrt{1-\frac{z+\log(1-z)}{\zc+\log(1-\zc)}}}
  \,.
\end{align}
This is now a root-finding problem for the function $g(\zc)$. We find  $\zc \approx -3.9052637703$.  Plugging this root back into \eqref{eq:eps_crit_z_min} yields
$\epsilon_{c} \approx 0.2159869288903$.

\subsection{A deeper look at the pitchfork bifurcation}

In retrospect, there were at least two reasons to expect that a pitchfork bifurcation could occur in our problem. The reasons have to do with the symmetry of the problem itself and the special structure of the associated vector field.  

First, recall that the  differential equation $\epsilon y'' = y y' - y$ and its boundary conditions $y(0)=1, y(1)=-1$ are left unchanged by the transformation $(x,y) \to (1-x,-y)$. If we apply this transformation a second time, we get back to $(x,y)$, as if toggling a switch or reflecting an image in a mirror twice. Such a transformation is known as a $\mathbb{Z}_2$ symmetry. Now it turns out that the occurrence of a pitchfork bifurcation in a boundary-value problem with $\mathbb{Z}_2$ symmetry is a codimension-1 phenomenon~\cite{mclachlan2018bifurcation}, which means that we should expect to see it in a generic one-parameter family of such problems. So we should not be surprised to find a pitchfork occurring in our problem as we vary its single parameter, $\epsilon$.

Second, recall from Sec.~\ref{sec:conserved} that the dynamical system~\eqref{eq:ODE_reduced} is Hamiltonian for $z < 1$. McLachlan and Offen~\cite{mclachlan2018bifurcation} give the generic bifurcations for Hamiltonian boundary-value problems and find that a pitchfork can be codimension-1 in planar problems even without the kind of symmetry seen in our system.

So if we had known what we know now, we should have expected a pitchfork all along. As always, everything becomes clearer in hindsight. Fortunately, before clarity comes, we have the pleasure of being surprised.

\section*{Acknowledgments}

This manuscript was written during the global COVID-19 pandemic, which has affected different members of society inequitably.  The authors were quite privileged to have had the opportunity to collaborate and continue to do research---opportunities not afforded to everyone.  At the same time, one impetus for this work was that S.H.S. taught his course on asymptotics and perturbation methods online, and L.C.S. became aware of Strogatz's lectures because of Twitter, and was able to follow along on YouTube (the problem at hand was discussed in~\cite{Strogatz-YT-NLBVP-2021}).  L.C.S. would therefore like to thank Twitter and YouTube for facilitating his initial involvement in this research.
We would also like to thank John Boyd, Mathieu Desroches, John Guckenheimer, Martin Krupa, Christian Kuehn, Ian Lizarraga, and Richard Rand for their helpful advice. 
In our numerical work we made use of the \texttt{python} package
\texttt{mpmath}~\cite{mpmath} for arbitrary precision calculation,
\texttt{numpy}~\cite{harris2020array},
\texttt{scipy}~\cite{2020SciPy-NMeth},
and
\texttt{matplotlib}~\cite{Hunter:2007}.

\appendix

\section{Higher-order asymptotics}
\label{sec:high-order-asympt}

Here we collect the long expressions obtained at higher orders of perturbation theory.
For the middle layer solution M, at first order, we find
\begin{align}
  \label{eq:Y1M}
  Y_{1}^{M} =
  \frac{1}{24} \sech^2\left(\frac{3 X}{4}\right) \Bigg[ & 16 \Li_2\left(-e^{-3 X/2}\right)
  +24 c_2 +3 X \left(3 X + 4+4 c_1 - 8 \log 2\right)
  \nonumber \\
  &{}
  + 8 \sinh \left(\frac{3 X}{2}\right) \left(2 \log \cosh \left(\frac{3 X}{4}\right)-1+c_1\right)
  \Bigg]
  \,,
\end{align}
where $X\equiv(x-1/2)/\epsilon$.
As discussed in the main text below \eqref{eq:Y1M_ODE}, after
matching, the integration constants take the values
$c_{1} = \pi^{2}/18$, and $c_{2} = 1+\log 4$.

For the B0 case, we condense the notation a bit by defining
$X\equiv x/\epsilon$ and
$\overline{X}\equiv X - \tanh^{-1} \frac{1}{2}$, so we can write
\begin{align}
  \label{eq:Y1B0}
  Y_{1}^{B0} =
  \frac{1}{4} \sech^2\overline{X} \Bigg[&
  4 c_2 + 2 \overline{X} \Big(\overline{X}
    + 1 + c_1 - \log 4\Big)
  \nonumber\\
  &{}+\sinh (2 \overline{X}) \Big(c_1 -1 +
  2 \log \cosh (\overline{X})\Big)
  + 2 \Li_2\left(-e^{-2 \overline{X}}\right)\Bigg]
  \,.
\end{align}
After matching we fix the integration constants as
$c_{1} = 1+\log 12$,
and
\begin{align}
  c_{2} = \frac{1}{8} \left(-4 \Li_2(-3)+(\log 3)^{2}+\frac{4}{3}\log 6912\right)
  \approx 2.594
  \,.
\end{align}
As discussed in the main text below \eqref{Mc_1}, the outer solution and the overlap solution cancel out when we form the composite solution, leaving only the inner solution. Thus the first-order composite solution for B0 is just
\begin{align}
  y_{c,1}^{B0} = -2 \tanh \overline{X} + \epsilon Y_{1}^{B0}
  \,,
\end{align}
where, as above, $\overline{X} = \frac{x}{\epsilon} -
\tanh^{-1}\frac{1}{2}$.

For the B1 solution, we can use the symmetry to write $y^{B1}(x) = - y^{B0}(1-x)$.

\bibliographystyle{siamplain}
\bibliography{references}

\end{document}